\documentclass[reqno]{amsart}

\usepackage[british]{babel}
\usepackage{graphicx}
\usepackage{amsfonts}
\usepackage{amssymb}
\usepackage{amsmath}
\usepackage[all]{xy}

\theoremstyle{plain}
\newtheorem{theorem}{Theorem}[section]
\newtheorem{lemma}[theorem]{Lemma}
\newtheorem{proposition}[theorem]{Proposition}
\newtheorem{corollary}[theorem]{Corollary}

\theoremstyle{definition}

\newtheorem{definition}[theorem]{Definition}

\newtheorem{remark}[theorem]{Remark}
\newtheorem{notation}[theorem]{Notation}
\newtheorem{example}[theorem]{Example}

\newcommand{\comp}{\raisebox{0.2mm}{\ensuremath{\scriptstyle{\circ}}}}
\newcommand{\defn}{\textbf}
\renewcommand{\mapsto}{\longmapsto}
\renewcommand{\to}{\longrightarrow}
\newcommand{\To}{\Longrightarrow}
\newcommand{\del}{\partial}
\newcommand{\iso}{\cong}

\renewcommand{\Im}{\ensuremath{\mathrm{Im}}}

\renewcommand{\ker}{\ensuremath{\mathrm{ker}}}
\newcommand{\Ker}{\ensuremath{\mathsf{Ker\,}}}

\newcommand{\A}{\ensuremath{\mathcal{A}}}
\newcommand{\B}{\ensuremath{\mathcal{B}}}
\newcommand{\C}{\ensuremath{\mathcal{C}}}
\newcommand{\D}{\ensuremath{\mathcal{D}}}
\newcommand{\E}{\ensuremath{\mathcal{E}}}
\newcommand{\Sc}{\ensuremath{\mathcal{S}}}
\newcommand{\Gp}{\ensuremath{\mathsf{Gp}}}
\newcommand{\Ab}{\ensuremath{\mathsf{Ab}}}
\newcommand{\AbGp}{\ensuremath{\mathsf{AbGp}}}
\DeclareMathOperator{\ab}{ab}
\DeclareMathOperator{\centr}{centr}
\DeclareMathOperator{\finab}{finab}
\DeclareMathOperator{\fincentr}{fincentr}

\newcommand{\XMod}{\ensuremath{\mathsf{XMod}}}
\newcommand{\SESeq}{\ensuremath{\mathsf{SESeq}}}
\newcommand{\SSeq}{\ensuremath{\mathsf{SSeq}}}

\newcommand{\Set}{\ensuremath{\mathsf{Set}}}
\newcommand{\Two}{\ensuremath{\mathsf{2}}}
\newcommand{\G}{\ensuremath{\mathbb{G}}}
\newcommand{\K}{\ensuremath{\mathbb{K}}}
\newcommand{\Z}{\ensuremath{\mathbb{Z}}}

\newcommand{\N}{\ensuremath{\mathbb{N}}}
\newcommand{\gnab}{\ensuremath{\text{<}}}
\newcommand{\Fin}{\ensuremath{\mathsf{Fin}}}
\newcommand{\End}{\ensuremath{\mathsf{End}}}
\newcommand{\Init}{\ensuremath{\widehat{\End}}}
\DeclareMathOperator{\cod}{cod}
\DeclareMathOperator{\dom}{dom}
\newcommand{\Ran}{\ensuremath{\mathsf{Ran}}}

\newcommand{\Lie}{\ensuremath{\mathsf{Lie}}}
\newcommand{\AbLie}{\ensuremath{\mathsf{AbLie}}}
\newcommand{\PXMod}{\ensuremath{\mathsf{PXMod}}}
\newcommand{\Fun}{\ensuremath{\mathsf{Fun}}}

\newcommand{\Ch}{\ensuremath{\mathsf{Ch}}}
\renewcommand{\hom}{\ensuremath{\mathrm{Hom}}}

\newcommand{\pr}{\ensuremath{\mathrm{pr}}}
\newcommand{\noproof}{\quad $\Box$}

\newcommand{\CExt}{\ensuremath{\mathsf{CExt}}}
\newcommand{\CExtnA}{\ensuremath{\mathsf{CExt}^{n}_{A}}}

\newcommand{\CExtk}{\ensuremath{\mathsf{CExt}^{k}_{\B}}}

\newcommand{\TExtk}{\ensuremath{\mathsf{TExt}^{k}_{\B}}}
\newcommand{\TExtkk}{\ensuremath{\mathsf{TExt}^{k+1}_{\B}\!}}
\newcommand{\Ext}{\ensuremath{\mathsf{Ext}}}

\newcommand{\Extk}{\ensuremath{\mathsf{Ext}^{k}\!}}
\newcommand{\Extkk}{\ensuremath{\mathsf{Ext}^{k+1}\!}}
\newcommand{\Extkn}{\ensuremath{\mathsf{Ext}^{k+n}\!}}

\newcommand{\ExtkkA}{\ensuremath{\mathsf{Ext}^{k+1}_{A}}\!}
\newcommand{\ExtknA}{\ensuremath{\mathsf{Ext}^{k+n}_{A}}\!}

\newcommand{\Arr}{\ensuremath{\mathsf{Arr}}}

\newcommand{\Arrk}{\ensuremath{\mathsf{Arr}^{k}\!}}
\newcommand{\Arrkk}{\ensuremath{\mathsf{Arr}^{k+1}\!}}

\newbox\skewpullbackbox
\setbox\skewpullbackbox=\hbox{\xy 0;<1mm,0mm>: \POS(4,0)\ar@{-} (0,4)
\ar@{-} (0,-4)
\endxy}
\def\skewpullback{\copy\skewpullbackbox}

\newbox\pullbackbox
\setbox\pullbackbox=\hbox{\xy 0;<1mm,0mm>: \POS(4,0)\ar@{-} (0,0) \ar@{-} (4,4)
\endxy}
\def\pullback{\copy\pullbackbox}

\newbox\pushoutbox
\setbox\pushoutbox=\hbox{\xy 0;<1mm,0mm>: \POS(0,4)\ar@{-} (0,0) \ar@{-} (4,4)
\endxy}

\newdir{>>}{{}*!/3.5pt/:(1,-.2)@^{>}*!/3.5pt/:(1,+.2)@_{>}*!/7pt/:(1,-.2)@^{>}*!/7pt/:(1,+.2)@_{>}}
\newdir{ >>}{{}*!/8pt/@{|}*!/3.5pt/:(1,-.2)@^{>}*!/3.5pt/:(1,+.2)@_{>}}
\newdir{ |>}{{}*!/-3.5pt/@{|}*!/-8pt/:(1,-.2)@^{>}*!/-8pt/:(1,+.2)@_{>}}
\newdir{ >}{{}*!/-8pt/@{>}}
\newdir{>}{{}*:(1,-.2)@^{>}*:(1,+.2)@_{>}}
\newdir{<}{{}*:(1,+.2)@^{<}*:(1,-.2)@_{<}}

\begin{document}

\title[On satellites in semi-abelian categories]{On satellites in semi-abelian categories:\\
Homology without projectives}

\author{Julia Goedecke}
\email{julia.goedecke@cantab.net}
\address{DPMMS\\
University of Cambridge\\
Wilberforce Road\\
Cambridge CB3 0WB\\
United Kingdom}

\author{Tim Van der Linden}
\email{tvdlinde@vub.ac.be}
\address{Vakgroep Wiskunde\\
Vrije Universiteit Brussel\\
Pleinlaan~2\\
1050~Brussel\\
Belgium}

\begin{abstract}
Working in a semi-abelian context, we use Janelidze's theory of generalised satellites to study universal properties of the Everaert long exact homology sequence. This results in a new definition of homology which does not depend on the existence of projective objects. We explore the relations with other notions of homology, and thus prove a version of the higher Hopf formulae. We also work out some examples.
\end{abstract}

\maketitle

\setcounter{tocdepth}{1}
\tableofcontents

\section{Introduction}

In his thesis \cite{Tomasthesis}, Everaert shows that, given the reflector $I\colon{\A\to \B}$ of a semi-abelian category~$\A$ with enough projectives to a Birkhoff subcategory $\B$ of~$\A$, any short exact sequence 
\[
\xymatrix{0 \ar[r] & K[f] \ar[r]^-{\Ker f} & B \ar[r]^-{f} & A \ar[r] & 0} 
\]
in~$\A$ induces a long exact sequence in $\B$,
\[
{\xymatrix@R=1ex{{\cdots} \ar[r] & H_{n+1}A \ar[r]^-{\delta^{n+1}_{f}} & K[H_{n}(f,I_{1})] \ar[r]^-{\gamma^{n}_{f}} & H_{n}B \ar[r]^-{H_{n}f} & H_{n}A \ar[r] & \cdots\\
{\cdots} \ar[r] & H_{2}A \ar[r]_-{\delta^{2}_{f}} & K[H_{1}(f,I_{1})] \ar[r]_-{\gamma^{1}_{f}} & H_{1}B \ar[r]_-{H_{1}f} & H_{1}A \ar[r] & 0,}}
\]
where the $H_{n}A=H_{n}(A,I)$ denote the homology of the object $A$ with coefficients in $I$, but $H_{n}(f,I_{1})$ is the homology of the extension $f$ with coefficients in $I_{1}$, the \emph{centralisation functor} associated with~$I$. This \defn{Everaert sequence}---a kind of generalised long Stallings-Stammbach sequence---no longer satisfies the classical abelian-categories properties of a long exact homology sequence. For instance, it is not functorial in the objects of the given short exact sequence: $K[H_{n}(f,I_{1})]$ need not be of the form $H_{n}K[f]$. 

We use Janelidze's theory of generalised satellites~\cite{Janelidze-Satellites} to arrive at a better understanding of this sequence's universal properties. Eventually this gives a way to compute homology using Kan extensions---as a limit---instead of basing it on higher Hopf formulae (as Everaert does) or simplicial resolutions (as, e.g., Barr and Beck do~\cite{Barr-Beck}). Thus we obtain a homology theory which also makes sense in a context where not enough projective objects are available. Our approach seems to be related to the work of Guitart and Van den Bril \cite{Guitart-Bril, Guitart:Anabelian} on homology using Kan extensions.

\subsection{Semi-abelian homology, Barr-Beck style}\label{Subsection-Barr-Beck-Definition}

In this paper, as in \cite{Tomasthesis, EverHopf, EGVdL, EverVdL2} and others, \emph{semi-abelian homology} studies the following classical situation. $\A$~is a semi-abelian category \cite{Janelidze-Marki-Tholen} (say, the category $\Gp$ of groups or $\Lie_{\K}$ of Lie algebras over a field $\K$ or $\PXMod$ of precrossed modules) and $\B$ a Birkhoff subcategory of $\A$ (the category $\Ab$ of abelian groups or $\AbLie_{\K}$ of abelian Lie algebras over $\K$ or $\XMod$ of crossed modules). Since the reflector $I\colon{\A\to\B}$ is not an exact functor, one is interested in its derived functors, as they capture some interesting homological information: integral homology of groups or homology of Lie algebras or of crossed modules.

A \defn{Birkhoff subcategory} $\B$ of a Barr-exact category $\A$ is a full reflective subcategory which is closed under subobjects and regular quotients \cite{Janelidze-Kelly}. For instance, a Birkhoff subcategory of a semi-abelian variety of universal algebras is the same as a subvariety. When $\A$ is a semi-abelian monadic category (e.g., a semi-abelian variety; see \cite{Gran-Rosicky:Monadic} for a precise characterisation), canonical regular-projective simplicial resolutions exist in $\A$, and we obtain the following Barr-Beck style \cite{Barr-Beck} notion of homology \cite{EverVdL2}: for any object $A$ of $\A$ and any $n\geq 0$,
\begin{equation}\label{Barr-Beck-Homology}
H_{n+1}(A,I)_\G=H_{n}NI\G A,
\end{equation}
where $I\colon{\A\to \B}$ is the reflector, $\G A$ is the simplicial resolution of $A$ obtained via the canonical forgetful/free comonad $\G$ on $\A$, and $N\colon{\Sc\B\to \Ch\B}$ is the \defn{Moore normalisation functor} which sends a simplicial object in $\B$ to its normalised chain complex. Note the dimension shift in~\eqref{Barr-Beck-Homology}; it is there for historical reasons: this is how, for example, homology of groups is numbered classically.

\subsection{Higher central extensions and the Hopf formulae}

It turns out that in the study of these homology objects, the concept of \emph{higher central extension} is fundamental. In \cite{EGVdL}, explicit Hopf formulae are proven which completely describe the $H_{n+1}(A,I)_\G$ in terms of centralisation of higher extensions. The most compact way to express their meaning seems to be that \emph{the $(n+1)$-st homology of $A$ measures the difference between the centralisation and the trivialisation of an $n$-fold presentation of~$A$.} Indeed, according to~\cite{EverHopf}, the Hopf formula of \cite[Theorem~8.1]{EGVdL} may be written as an isomorphism
\[
H_{n+1}(A,I)_{\G}\cong K^{n+1}[I_{n}p\to T_{n}p]
\]
where $p$ is an $n$-fold presentation of $A$ and $n\geq 1$. The notions of central and trivial extension and the meaning of all ingredients of this formula will be explained in Section~\ref{Section-Preliminaries}.

\subsection{The Hopf formulae as a definition of homology}\label{Subsection-Everaert-Definition} 

This idea---to explain homology objects in terms of higher-dimensional central and trivial extensions---is further pursued by Everaert in \cite{Tomasthesis} and \cite{EverHopf}, where he works out a new notion of homology based on the right hand side of the Hopf formula isomorphism: there \emph{by definition}, $H_{n+1}(A,I)=K^{n+1}[I_{n}p\to T_{n}p]$, for $n\geq 1$ and any $n$-fold presentation $p$ of $A$. Note how the comonad $\G$ is dropped from the notation. In fact, as explained in \cite{EverHopf}, this approach, using higher presentations of an object, is much closer to Hopf's original insights than the use of simplicial resolutions. When the underlying category~$\A$ is semi-abelian and monadic, the higher Hopf formulae become $H_{n+1}(A,I)_\G\cong H_{n+1}(A,I)$, the equivalence between the two notions of homology. But Everaert's theory works as soon as $\A$ is semi-abelian with enough projectives, while it is still powerful enough to obtain interesting results: no monadicity condition on $\A$ is needed to obtain, say, a long exact homology sequence.

\subsection{A third approach: homology via satellites}

It turns out that the universal properties of the Everaert homology sequence completely determine an underlying homology theory, and these universal properties may be taken as a new definition of homology. The advantage of such an approach is that the existence of projective objects is no longer fundamental, and new homological techniques are obtained. This is the subject of the present paper.

Our theory is based on Janelidze's general notion of \emph{satellites} \cite{Janelidze-Satellites}, which give a way to compute homology objects step by step: the $(n+1)$-st homology $H_{n+1}$ is obtained out of $H_{n}$ as a Kan extension. This makes it possible to define homology using limits alone. But when the surrounding category has enough projectives, the resulting notion is still equivalent to Everaert's---an equivalence which may be interpreted as a version of the higher Hopf formulae valid in this context.

\subsection{Structure of the text}

In Section~\ref{Section-Preliminaries} we briefly sketch some of the basic definitions and properties used throughout the text. Section~\ref{Section-Satellites} is devoted to the definition of satellites and the proof that the homology objects in the sense of~\ref{Subsection-Everaert-Definition} (and hence also in the sense of \ref{Subsection-Barr-Beck-Definition}) are satellites. The main results here are Proposition~\ref{Proposition-Satellites} (which gives $H_{n+1}(-,I)$ as a satellite of $H_{n}(-,I_{1})$) and Theorem~\ref{Theorem-Satellites} (which gives $H_{n+1}(-,I)$ as a satellite of $I_{n}$). In Section~\ref{Section-Without-Projectives} satellites are used to define homology. In Section~\ref{Section-Homology-With-Projectives} the consequences of this definition are explored in the situation where enough projective objects do exist. In that case, homology can be calculated in a new way, as the limit of a certain small diagram involving a projective presentation.

\section{Preliminaries}\label{Section-Preliminaries}

\subsection{Semi-abelian categories}

First of all, we shall not limit ourselves to \defn{semi-abelian} categories (which are pointed, Barr exact and Bourn protomodular with binary coproducts \cite{Janelidze-Marki-Tholen, Borceux-Bourn}) but choose \emph{pointed exact protomodular} categories as the basic context. All constructions we borrow from \cite{EGVdL} and \cite{Tomasthesis} and which take place in a semi-abelian category still work in pointed exact protomodular ones---though they need not have coproducts, these categories still have cokernels of kernels (see \cite[Corollary~4.1.3]{Borceux-Bourn}). Since the rest of our theory also does not need coproducts, it seems unnecessary to require their existence.

\subsection{Higher-dimensional arrows}

We are interested in the chain of inclusions of full subcategories
\[
\Arrk\A \supset \Extk\A \supset \CExtk\A \supset \TExtk\A
\]
where $k\geq 1$, $\A$ is a pointed exact protomodular category and $\B$ a Birkhoff subcategory of~$\A$. The category $\Arrk\A$ consists of \defn{$k$-dimensional arrows} in $\A$: $\Arr^{0}\!\A=\A$, $\Arr^{1}\!\A=\Arr\A$ is the category of arrows $\Fun(\Two,\A)$ where $\Two$ is generated by a single map ${\varnothing\to \{\varnothing\}}$, and $\Arrkk\A=\Arr\Arrk\A$. Thus a~\defn{double arrow} is a commutative square in $\A$, a $3$-arrow is a commutative cube, and a $k$-arrow is a commutative $k$-cube. Clearly, $\Arrk\A$ is also pointed exact protomodular. The functor $\ker\colon{\Arrkk\A\to \Arrk\A}$ maps a $(k+1)$-arrow $a$ to its kernel $K[a]$, and a morphism $(f',f)$ between $(k+1)$-arrows $b$ and $a$ to the induced morphism between their kernels.
\[
\xymatrix{K[b] \ar[d]_-{\ker (f',f)} \ar@{ |>->}[r]^-{\Ker b} & B' \ar[d]_-{f'} \ar[r]^-{b} \ar@{}[dr]|{\Downarrow} & B \ar[d]^-f\\
K[a] \ar@{ |>->}[r]_-{\Ker a} & A' \ar[r]_-a & A}
\]
Repeating it $n$ times gives a functor $\ker^{n}\colon{\Arr^{k+n}\!\A\to \Arrk\A}$ which sends a $(k+n)$-arrow $a$ to the object $K^{n}[a]$ of $\Arrk\A$.

\subsection{Extensions}

A \defn{$0$-extension in~$\A$} is an object of $\A$ and a \defn{$1$-extension} is a regular epimorphism in $\A$. For $k\geq 2$, a \defn{$k$-extension} is an object $(f',f)$ of $\Arrk\A$ such that all arrows in the induced diagram
\begin{equation}\label{Diagram-Extension}
\vcenter{\xymatrix{B' \ar@/^/@{-{ >>}}[drr]^-{f'} \ar@/_/@{-{ >>}}[drd]_-{b} \ar@{. >>}[rd]|r \\
& P \ar@{}[rd]|<{\pullback} \ar@{.{ >>}}[r] \ar@{.{ >>}}[d] & A' \ar@{-{ >>}}[d]^-a \\
& B \ar@{-{ >>}}[r]_-f & A}}
\end{equation}
are $(k-1)$-extensions. Here $P$ is the pullback of $a$ and $f$. The $k$-extensions determine a full subcategory $\Extk\A$ of $\Arrk\A$. A $2$-extension is better known as a \defn{double extension}, and $\Ext\A=\Ext^{1}\!\A$. When we say that a sequence is exact in $\Extk\A$, we mean that it is an exact sequence in $\Arrk\A$, and the objects are $k$-extensions. Given a short exact sequence
\begin{equation}\label{Equation-Short-Exact-Sequence}
\xymatrix@1{0 \ar[r] & K[f] \ar@{ |>->}[r]^-{\Ker f} & B \ar@{ >>}[r]^-{f} & A \ar[r] & 0} 
\end{equation}
in $\Arrk\A$, all three objects are $k$-extensions if and only if the map $f$ is a $(k+1)$-extension, by~\cite[Proposition~3.9]{EGVdL}.

Roughly, the idea behind this definition of $k$-extensions is the following: suppose we are given a double extension $(f',f)$ of an object $A$ of $\A$ as in Diagram~\eqref{Diagram-Extension}, and let $\alpha$ be any element of $A$. Then in addition to the existence of elements $\beta$ of $B$ and $\alpha'$ of $A'$ such that $f(\beta)=\alpha$ and $a(\alpha')=\alpha$, there is also an element $\beta'\in B'$ such that $b(\beta')=\beta$ and $f'(\beta')=\alpha'$, whichever $\beta$ and $\alpha'$ were chosen.

\subsection{The Galois structures $\Gamma_{k}$}

A Birkhoff subcategory $\B$ of a pointed exact protomodular category $\A$ together with its reflector $I\colon{\A\to \B}$ and the classes of extensions in $\A$ and $\B$ forms a \defn{Galois structure} in the sense of Janelidze~\cite{Janelidze:Precategories}. With respect to this Galois structure $\Gamma_{0}$, there is a notion of \defn{central extension} such that the full subcategory $\CExt_{\B}\A$ of $\Ext\A$ determined by the central extensions is again reflective. Its reflector $I_{1}\colon{\Ext\A\to \CExt_{\B}\A}$, together with the classes of extensions in $\Ext\A$ and in $\CExt_{\B}\A$ (i.e., double extensions in $\A$, and double extensions with central domain and codomain), in turn determines a Galois structure $\Gamma_{1}$. Inductively, this defines a family of Galois structures $(\Gamma_{k})_{k\in\N}$, each of which gives rise to a notion of central extension which determines the next structure. (See Subsection~\ref{Subsection-Reflectors} and~\cite{Tomasthesis, EverHopf, EGVdL} for more details.) In particular, for every $k\geq 1$ we obtain a reflector
\[
I_{k}\colon{\Extk\A\to\CExtk\A,}
\]
left adjoint to the inclusion $\CExtk\A\subset\Extk\A$.

\subsection{The reflectors $I_{k}$}\label{Subsection-Reflectors}

We will not spend too much time in this paper explaining the Galois structures $\Gamma_k$ in detail, but only sketch the construction of the reflectors $I_{k}\colon{\Extk\A\to\CExtk\A}$. We can view the reflector $I=I_0$ as a functor $I\colon {\A \to \A}$. Let $\eta\colon{1_{\A}\To I}$ be the unit of the adjunction associated with $I$. Then we have another functor $J\colon {\A \to \A}$, given by $JA=K[\eta_A]$, which fits into the following short exact sequence of functors.
\[
 \xymatrix{0 \ar[r] & J \ar@{{ |>}->}[r]^{\mu} & 1_{\A} \ar@{ >>}[r]^{\eta} & I \ar[r] & 0}
\]
 From this, we build a similar short exact sequence of functors $\Ext\A \to \Ext\A$ as follows. (The construction is made pointwise in $\Arr\A$, which has good categorical properties, but the result turns out to be an extension.) Consider an extension $f\colon {B\to A}$ and its kernel pair $(\pi_1, \pi_2)$. Write $J_1[f]=K[J\pi_1]$ and $J_1f\colon {J_1[f]\to 0}$. 
\[
\xymatrix{J_1[f]=K[J\pi_1] \ar@{ |>->}[d] \ar@{ |>->}[r]^-{\Ker J\pi_1} \ar@{}[dr]|<<{\pullback} & JR[f] \ar@{ |>->}[d]_-{\mu_{R[f]}} \ar@<.5ex>[r]^-{J\pi_1} \ar@<-.5ex>[r]_-{J\pi_2} & JB \ar@{ |>->}[d]^-{\mu_B}\\
K[f]=K[\pi_1] \ar@{ |>->}[r]_-{\Ker \pi_1} & R[f] \ar@<.5ex>[r]^-{\pi_1} \ar@<-.5ex>[r]_-{\pi_2} & B}
\]
This clearly determines a functor $J_1\colon{ \Ext\A \to \Ext\A}$. Note that $\pi_2 \comp \Ker \pi_1=\Ker f$, and the left hand square is a pullback. We define the map $\mu^1_f\colon {J_1f\to f}$ as in the left hand square below.
\[
\xymatrix{J_1[f] \ar[rr] \ar@<1ex>@{}[rr]^-{\mu_B\comp J\pi_2\comp \Ker J\pi_1} \ar@{ >>}[d]_-{J_1f}  \ar@{}[drr]|{\stackrel{\mu^1_f}{\To}} && B \ar@{ >>}[d]^-{f}\\
0 \ar[rr] & &A}
\qquad\qquad
\xymatrix{B \ar@{ >>}[r]^-{\rho^1_f} \ar@{ >>}[d]_-{f}\ar@{}[dr]|{\stackrel{\eta^1_f}{\To}} & I_1[f] \ar@{ >>}[d]^-{I_1f}\\
A \ar@{=}[r] & A}
\]
Note that the composition $\mu_B\comp J\pi_2\comp\Ker J\pi_1$ is a normal monomorphism, so we can take cokernels, yielding the right hand square. Since $\mu^{1}_{f}$ is the kernel of its cokernel, we obtain the short exact sequence
\[
\xymatrix{0 \ar[r] & J_1 \ar@{ |>->}[r]^{\mu^1} & 1_{\Ext\A} \ar@{ >>}[r]^{\eta^1} & I_1 \ar[r] &0}
\]
of functors $\Ext\A \to \Ext\A$. This process may be repeated inductively to obtain the functors $J_k$ and $I_k$ from $\Extk\A$ to $\Extk\A$. For $k\geq 1$ and a $k$-extension $f$, we often call the extension $I_kf$ the \defn{centralisation of $f$}.

\begin{remark}\label{Remark-Bang-Reflector}
Given a $k$-extension $A$, for $k\geq 0$, the centralisation of the $(k+1)$-extension $!_A\colon {A\to 0}$ turns out to be $I_{k+1}!_A\colon {I_kA\to 0}$.
\end{remark}

The following is also often useful, and quite easy to show using the $3\times 3$-Lemma and the strong (extension)-Birkhoff property \cite[Definition~2.5]{EGVdL} satisfied by the category of central $k$-extensions (see \cite[Lemma~4.3]{EGVdL}).

\begin{lemma}\label{Lemma-II_1=I}\cite[Lemma~6.2]{EGVdL} For a $(k+1)$-extension $f\colon {B\to A}$, we have
\[
I_kI_{k+1}f= I_kf\colon {I_kB\to I_kA},
\]
i.e., $I_k(I_{k+1}[f])=I_kB$. \noproof
\end{lemma}

\begin{remark}\label{Remark-Sparse-Jn}
Given any $k$-extension $f$, the only object of $J_kf$ which is non-zero is $\dom^k\!J_kf$, the ``initial'' object of the $k$-cube $J_kf$. This follows easily from the inductive construction of $J_kf$. Thus we have $\dom^k\!J_kf=K^k[J_kf]$ for any $k$-extension $f$. 
\end{remark}

\subsection{Trivial extensions}

A \emph{trivial extension} is a special kind of central extension: a $(k+1)$-extension $f\colon{B\to A}$ is \defn{trivial} (with respect to the Galois structure $\Gamma_{k}$) when it is the pullback of its reflection $I_{k}f\colon{I_{k}B\to I_{k}A}$ into $\CExtk\A$ along the unit $\eta^{k}_{A}\colon{A\to I_{k}A}$ at $A$ of the reflector~$I_{k}$. The trivial $(k+1)$-extensions of $\A$ form a reflective subcategory $\TExtkk\A$ of $\Extkk\A$; the reflector
\[
T_{k+1}\colon{\Extkk\A\to \TExtkk\A}
\]
maps an extension $f$ to the pullback $T_{k+1}f\colon{T_{k+1}[f]\to A}$ of $I_{k}f$ along $\eta^{k}_{A}$, the \defn{trivialisation of~$f$}.
\begin{equation}\label{Diagram-Compare-I-T}
\vcenter{\xymatrix@!0@=3,5em{&&&&& A \ar@{ >>}[rd]^{\eta^k_A}\\
B \ar@/^1.5pc/@{ >>}[rrrrru]^-{f} \ar@/_1.5pc/@{ >>}[rrrrrd]_-{\eta^k_B}
\ar@{ >>}[rr]|-{\rho^k_f} && I_{k+1}[f] \ar@/^0.5pc/@{ >>}[rrru]|(.4){I_{k
+1}f} \ar@/_0.5pc/@{ >>}[drrr]|(.36){\eta^k_{I_{k+1}[f]}} \ar@{.{ >>}}[rr] &&
T_{k+1}[f]  \ar@{}[r]|<{\skewpullback} \ar@{ >>}[ru]|{T_{k+1}f}
\ar@{ >>}[rd] && I_kA\\
&&&&& I_kB \ar@{ >>}[ru]_-{I_kf}& }}
\end{equation}
Thus we obtain a comparison map $r^{k+1}_f\colon {I_{k+1}[f] \to T_{k+1}[f]}$, which is a $(k+1)$-extension by the strong (extension)-Birkhoff property \cite[Definition~2.5]{EGVdL} of the reflector $I_{k}$ and \cite[Lemma~3.8]{EGVdL}. This gives a $(k+2)$-extension ${I_{k+1}f \to T_{k+1}f}$.

\begin{remark}\label{Remark-Central-Extensions}
 The Galois-theoretic definition of a central extension \cite{Janelidze:Precategories} says that an extension $f\colon {B\to A}$ is central if and only if there is an extension $g\colon {\overline A\to A}$ such that the pullback $\overline{f}\colon {\overline B\to \overline A}$ of $f$ along $g$ is a trivial extension.
\end{remark}

\subsection{Projective presentations}\label{Subsection-Projective-Presentations}

An object of $\Arrk\A$ is \defn{extension-projective} if it is projective with respect to the class of $(k+1)$-extensions. A $(k+1)$-extension $f\colon {B\to A}$ is called a \defn{(projective) presentation} of~$A$ when the object $B$ is extension-projective. A $(k+n)$-extension $f\colon {B\to A}$ is called an \defn{$n$-fold presentation}, or just \defn{$n$-presentation}, when the object $B$ is extension-projective and $A$ is an $(n-1)$-presentation. (A $1$-presentation is just a projective presentation as above.) Given an object $A$ of~$\Extk\A$, a \defn{$n$-fold presentation $p$ of~$A$} is an $n$-fold presentation with $\cod^{n}\!p =A$---the ``terminal object'' of the $n$-cube $p$ in $\Extk\A$ is~$A$.

\subsection{Key results}

Now we have provided definitions for all elements of the Hopf formula---the isomorphism
\[
H_{n+1}(A,I)_{\G}\cong K^{n+1}[I_{n}p\to T_{n}p],
\]
valid for any $n$-fold presentation $p$ of $A$ and any $n\geq 1$~\cite{EverHopf, EGVdL}. The crucial point here is that the information in the higher homology objects is entirely contained in higher-dimensional versions $I_k\colon{\Extk\A\to \CExtk\A}$ of the reflector $I\colon {\A\to\B}$. In this section and in Section~\ref{Section-Satellites}, we use homology \emph{defined} via the Hopf formulae, as in Section~\ref{Subsection-Everaert-Definition}: for any $k$-extension $A$ and an $n$-fold presentation $p$ of $A$, we define
\[
 H_{n+1}(A,I_k)=K^{n+1}[I_{k+n}p\to T_{k+n}p].
\]

\begin{remark}\label{Remark-Hopf-Formula-Different-Version}
Notice that in \cite{Tomasthesis,EGVdL}, the Hopf formula has the form
\[
\frac{J_kP_n\cap K^n[p]}{K^n[J_{n+k}p]},
\]
where $P_n$ is the ``initial'' $k$-extension in the $n$-cube representing $p$, and $K^n[p]=\bigcap_{i=0}^nK[p_i]$ is the intersection of all maps $p_i$ with domain $P_n$ in $p$. Everaert shows in \cite[Remark~5.12]{EverHopf} that this is indeed equivalent to the form we are using. 
\end{remark}

\begin{theorem}\cite[Theorem~2.4.2]{Tomasthesis}\label{Theorem-Everaert-Sequence}
For any $k\geq 0$, any short exact sequence~\eqref{Equation-Short-Exact-Sequence} in $\Extk\A$ induces a long exact homology sequence
\begin{equation}\label{Everaert-Sequence}
\resizebox{\textwidth}{!}
{\xymatrix@R=1ex{{\cdots} \ar[r] & H_{n+1}(A,I_{k}) \ar[r]^-{\delta^{n+1}_{f}} & K[H_{n}(f,I_{k+1})] \ar[r]^-{\gamma^{n}_{f}} & H_{n}(B,I_{k}) \ar[r]^-{H_{n}(f,I_{k})} & H_{n}(A,I_{k}) \ar[r] & \cdots\\
{\cdots} \ar[r] & H_{2}(A,I_{k}) \ar[r]_-{\delta^{2}_{f}} & K[H_{1}(f,I_{k+1})] \ar[r]_-{\gamma^{1}_{f}} & H_{1}(B,I_{k}) \ar[r]_-{H_{1}(f,I_{k})} & H_{1}(A,I_{k}) \ar[r] & 0}}
\end{equation}
in $\Extk\A$.
\end{theorem}
\begin{proof}
A proof of this theorem in its full generality is given in \cite{Tomasthesis}. However, when we restrict ourselves to the monadic case it becomes relatively easy to understand why the sequence takes this shape. So suppose that $\A$ is a semi-abelian monadic category and $\G$ the induced comonad on $\Extk\A$. This comonad produces canonical simplicial resolutions $\G A$ and $\G B$ of $A$ and $B$ and, by functoriality, also a simplicial resolution $\G f$ of $f$. The Everaert Sequence~\eqref{Everaert-Sequence} is the long exact homology sequence (see \cite[Corollary~5.7]{EverVdL2}) obtained from the short exact sequence of simplicial objects
\[
\xymatrix{0 \ar[r] & K[I_k\G f] \ar@{ |>->}[r] & I_k\G B \ar@{ >>}[r]^-{I_k\G f} & I_k\G A \ar[r] & 0;} 
\]
it remains to be shown that $H_{n-1}K[I_{k}\G f] = K[H_n(f,I_{k+1})]$ for all $n\geq 1$. (Remember the dimension shift in Equation~\eqref{Barr-Beck-Homology}.) Now degree-wise, the $(k+1)$-extension 
\[
I_{k+1}\G f\colon {I_{k+1}[\G f]\to \G A}
 \]
 is a split epimorphic central extension: it is a centralisation, and $\G A$ is degree-wise projective. Via~\cite[Proposition~4.5]{EGVdL}, this implies that, degree-wise, it is a trivial extension. This means that $I_{k+1}\G f$ is the pullback of $I_k\G f$ along the unit $\eta_{\G A}^{k}\colon {\G A \to I_k\G A}$, which in turn implies that $K[I_k\G f]$ is the kernel $K[I_{k+1}\G f]$ of $I_{k+1}\G f$. Since, $\G A$ being a simplicial resolution, $H_n\G A=0$ for all $n\geq 1$, the long exact homology sequence induced by the short exact sequence of simplicial objects
\[
\xymatrix{0 \ar[r] & K[I_{k+1}\G f] \ar@{ |>->}[r] & I_{k+1}[\G f] \ar@{ >>}[r]^-{I_{k+1}\G f} & \G A \ar[r] & 0} 
\]
gives the needed isomorphism $H_{n-1}K[I_{k+1}\G f] \cong K[H_n(f,I_{k+1})]$.
\end{proof}

Note that in \cite{Tomasthesis}, this sequence has a slightly different appearance: there it contains the objects $\dom H_{n}(f,I_{k+1})$ instead of $K[H_{n}(f,I_{k+1})]$ for $n\geq 2$. But the codomain of $H_{n}(f,I_{k+1})$ is zero (because a $k$-extension $J_kf$ is only non-zero in the very top corner of the $k$-cube representing $J_kf$, hence $I_k$ only changes the very top object of a $k$-extension), so its domain coincides with its kernel. For us, the sequence in its present, more uniform, shape will be easier to work with.

\begin{corollary}(cf.\ \cite[Theorem~6.4]{EGVdL})\label{Corollary-Homology-of-extension}
For any $n\geq 2$, $k\geq 0$ and any projective presentation $p\colon {P\to A}$ of a $k$-extension $A$,
\[
K[H_{n} (p,I_{k+1})]\cong H_{n+1} (A,I_{k}).
\]
\end{corollary}
\begin{proof}
It suffices to note that in the Everaert Sequence~\eqref{Everaert-Sequence}, all $H_{n+1} (P,I_{k})$ are zero, because $P$ is projective.
\end{proof}

This shows how the degree of the homology may be lowered from $n+1$ to $n$ by raising the degree of the reflector from $k$ to $k+1$.

\section{Satellites and homology}\label{Section-Satellites}

This section gives an analysis of homology in terms of satellites. Again we mean homology as defined in Section~\ref{Subsection-Everaert-Definition}. We start by stating the main definitions. Then, in Subsection~\ref{Subsection-H_{n+1}}, we interpret $H_{n+1}(-,I_k)$ (together with the connecting map~$\delta^{n+1}$) as a satellite of $H_{n}(-,I_{k+1})$. In Subsection~\ref{Subsection-Main-Theorem} we prove the main theorem of this section: a formula which gives $H_{n+1}$ in terms of~$I_{n}$. Finally in Subsection~\ref{Subsection-Gamma} we explain how the situation is entirely symmetric, in that the connecting map $\gamma^{n}$ also arises as a pointwise satellite.

\subsection{Satellites and pointwise satellites}\label{Subsection-Definitions}

Modulo a minor terminological change, the following definition is due to Janelidze.

\begin{definition}\label{Definition-Satellites}\cite[Definition~2]{Janelidze-Satellites}
Let $I'\colon{\A'\to\B'}$ be a functor. A \defn{left satellite $(H,\delta)$ of $I'$ (relative to $F\colon{\A'\to\A}$ and $G\colon{\B'\to\B}$)} is a functor $H\colon{\A\to \B}$ together with a natural transformation $\delta\colon {HF\To GI'}$
\[
\xymatrix@!0@=3,5em{& \A' \ar[ld]_-{F} \ar[rd]^-{I'}\\
\A \ar@{.>}[rd]_-{H} && \B' \ar[ld]^-{G}\\
& \B \ar@{:>}^-{\delta} "2,1";"2,3"}
\]
universal amongst such, i.e., if there is another functor $L\colon{\A\to\B}$ with a natural transformation $\lambda\colon{LF\To GI'}$, then there is a unique natural transformation $\mu\colon{L\To H}$ satisfying $\delta\comp \mu_{F}=\lambda$. This means that $(H,\delta)$ is the right Kan extension $\Ran_{F}GI'$ of the functor $GI'$ along $F\colon{\A'\to \A}$.
\end{definition}

This makes it possible to compute derived functors in quite diverse situations. The following example, borrowed from \cite{Janelidze-Satellites}, explains how satellites may be used to capture homology in the classical abelian case.

\begin{example}\label{Example-Abelian}
In the abelian context, the $(n+1)$-st homology functor $H_{n+1}$ may be seen as a left satellite of $H_n$. For instance, let $\A=\B'$ and $\B$ be categories of modules and $G\colon{\A\to \B}$ an additive functor. Then $G=H_0(-,G)$. Let $\SESeq\A$ be the category of short exact sequences
\[
\xymatrix{0 \ar[r] & K \ar@{ |>->}[r]^-{k} & B \ar@{ >>}[r]^-{f} & A \ar[r] & 0} 
\]
in $\A$, the functor $I'\colon{\SESeq \A\to \A}$ the projection $\pr_1$ that maps a sequence $(k,f)$ to the object~$K$, and $F\colon{\SESeq \A\to \A}$ the projection $\pr_3$ that maps $(k,f)$ to $A$. Let $H\colon{\A\to \B}$ be the first homology functor $H_1(-,G)$. We obtain a satellite diagram
\[
\xymatrix@!0@=3,5em{& {\SESeq \A} \ar[ld]_-{\pr_3} \ar[rd]^-{\pr_1}\\
\A \ar@{.>}[rd]_-{H_1(-,G)} && \A \ar[ld]^-{H_0(-,G)}\\
& \B \ar@{:>}^-{\delta} "2,1";"2,3"}
\]
where the natural transformation $\delta=(\delta_{(k,f)})_{(k,f)\in|\SESeq\A|}$ consists of the connecting maps from the (classical) long exact homology sequence 
\[
\resizebox{\textwidth}{!}
{\xymatrix{{\cdots} \ar[r] &  H_1K \ar[r]^-{H_1k} & H_1B \ar[r]^-{H_1f} & H_1A \ar[r]^-{\delta_{(k,f)}} &  H_0K \ar[r]^-{H_0k} & H_0B \ar[r]^-{H_0f} & H_0A \ar[r] & 0.}}
\]
The universality of the Kan extension follows from the universality of the long exact homology sequence amongst similar sequences and may for instance be shown as follows. Given any functor $L\colon{\A\to \B}$ and any natural transformation
\[
\lambda\colon{L\comp \pr_3\To H_0(-,G)\comp \pr_1},
\]
we will construct the component at an object $A\in|\A|$ of the needed natural transformation
\[
{L\To H_1(-,G)}
\]
by using a projective presentation $p\colon {P\to A}$ of $A$. Let $k\colon{K\to P}$ be the kernel of this projective presentation of $A$. Since $H_1P$ is zero (as $P$ is projective), the exactness of the long homology sequence induced by $(k,p)$ says that $\delta_{(k,p)}\colon {H_1A\to H_0K}$ is the kernel of $H_0k$. Then the string of equalities
\[
H_0k\comp \lambda_{(k,p)}\overset{(1)}{=}\lambda_{(1_P,!_P)}\comp L!_A\overset{(2)}{=}H_0(\gnab_P)\comp \lambda_{(1_0,1_0)}\comp L!_A\overset{(3)}{=}0
\]
yields the needed factorisation ${LA\to H_1A}$: $(1)$ expresses the naturality of $\lambda$ at the upper, downward-pointing morphism of the diagram
\begin{equation}\label{Subdiagram-Exact-Sequences}
\vcenter{\xymatrix{0 \ar[r] & K \ar@{ |>->}[r]^-{k} \ar[d]_-{k} \ar@{}[dr]|{\Downarrow} & P \ar@{=}[d]^{1_P} \ar@{ >>}[r]^-{p} \ar@{}[dr]|{\Downarrow} & A \ar[d]^{!_A} \ar[r] & 0\\
0 \ar[r] & P \ar@{=}[r]_-{1_P} \ar@{}[dr]|{\Uparrow}& P \ar@{ >>}[r]_-{!_P}\ar@{}[dr]|{\Uparrow}  & 0 \ar[r] & 0\\
0 \ar[r] & 0 \ar@{=}[r]_-{1_0} \ar[u]^-{\gnab_P} & 0 \ar[u]_{\gnab_P} \ar@{=}[r]_-{1_0} & 0 \ar[u]_{1_0} \ar[r] & 0}}
\end{equation}
in $\SESeq\A$, while $(2)$ follows from $\lambda_{(1_P,!_P)} = \lambda_{(1_P,!_P)} \comp L1_{0} = H_0(\gnab_P) \comp \lambda_{(1_0,1_0)}$, which is the naturality of $\lambda$ at the lower, upward-pointing morphism; the last equality $(3)$ holds because ${H_00=0}$.

Note that, as such, this example does \emph{not} follow the terminology of Definition~\ref{Definition-Satellites}. From its point of view one is tempted to call $H$ a left satellite of $G$ (rather than a satellite of $I'$), and actually this is how the definition appears in the paper \cite{Janelidze-Satellites}. But the situation we shall be considering in this paper demands the change in terminology, and the present example may easily be modified to comply with Definition~\ref{Definition-Satellites}.

Indeed, the functor $G$ may be lifted to a functor
\[
\SSeq  H_0(-,G)\colon{\SESeq\A\to \SSeq\B}
\]
where the latter category consists of short (not necessarily exact) sequences in~$\B$. Together with the obvious projection $\pr_{1}\colon {\SSeq\B\to \B}$ (such that $H_0(-,G)\comp \pr_1=\pr_1\comp \SSeq H_0(-,G))$, this gives us the satellite diagram
\[
{\xymatrix@!0@=3,5em{& {\SESeq \A} \ar[ld]_-{\pr_3} \ar[rd]^-{\SSeq  H_0(-,G)}\\
\A \ar@{.>}[rd]_-{H_1(-,G)} && \SSeq\B. \ar[ld]^-{\pr_1}\\
& \B \ar@{:>}^-{\delta} "2,1";"2,3"}}
\]
Whereas such a viewpoint may seem rather far-fetched in the abelian case, it is the only one still available when the context is widened to semi-abelian categories.
\end{example}

In practice, satellites may almost always be computed explicitly using limits---namely, as pointwise Kan extensions. Then the definition given above is strengthened as follows.

\begin{notation}
 Let $A$ be an object of $\A$. We denote by $(A\downarrow F)$ the category of elements of the functor $\hom(A,F-)\colon{\A'\to \Set}$: its objects are pairs $(A',\alpha\colon{A\to FA'})$, where $A'$ is an object of $\A'$ and $\alpha$ is a morphism in $\A$, and its morphisms are defined in the obvious way (cf.\ \cite[Theorem 3.7.2]{Borceux:Cats}). The forgetful functor $U\colon{(A\downarrow F)\to \A'}$ maps a pair $(A',\alpha)$ to $A'$. The natural transformation $(H,\delta)$ now induces a cone $\overline{\delta}$ on $G I' U\colon{(A\downarrow F)\to \B}$ with vertex $HA$ defined by
\[
\overline{\delta}_{(A',\alpha\colon{A\to FA'})}=\delta_{A'}\comp H\alpha\colon\xymatrix@1{HA \ar[r]^-{H\alpha} & HFA' \ar[r]^-{\delta_{A'}} & GI'A'=GI'U(A',\alpha).}
\]
\end{notation}

\begin{definition}\label{Definition-Pointwise-Satellite}
A left satellite $(H,\delta)$ of $I'$ relative to $F\colon{\A'\to\A}$ and $G\colon{\B'\to\B}$ is called \defn{pointwise} when it is pointwise as a Kan extension, i.e., for every object $A$ of $\A$, the cone $(HA,\overline{\delta})$ on $G I' U\colon{(A\downarrow F)\to \B}$ is a limit cone.
\end{definition}

To check that a pair $(H,\delta)$ is a pointwise satellite it is not necessary to prove its universality as in Definition~\ref{Definition-Satellites}, but it suffices to check the limit condition from Definition~\ref{Definition-Pointwise-Satellite}; see, for example, Mac\,Lane \cite[Theorem X.3.1]{MacLane}.

\subsection{$H_{n+1}(-,I_{k})$ as a satellite of $H_{n}(-,I_{k+1})$}\label{Subsection-H_{n+1}}

We are now ready to prove the first main result of this paper: we focus on the universal properties of the Everaert Sequence~\eqref{Everaert-Sequence}, and prove that they allow us to interpret the $(n+1)$-st homology with coefficients in $I_k$ as a satellite of the $n$-th homology with coefficients in $I_{k+1}$.

\begin{lemma}\label{Lemma-iota}\label{Lemma-iota-plus}
For $n\geq 1$, $k\geq 0$ and $A\in |\Extk\A|$,
\[
K[H_{n}(!_{A}\colon{A\to 0},I_{k+1})]=H_{n}(A,I_{k}).
\]
\end{lemma}
\begin{proof}
This follows from the exactness of the Everaert Sequence~\eqref{Everaert-Sequence} and the fact that all $H_n(0,I_k)$ are zero.
\end{proof}

\begin{lemma}\label{Lemma-Gamma}
For all $n\geq 1$, $k\geq 0$ and $f\colon{B\to A}\in |\Extkk\A|$,
\[
\gamma^{n}_{f} = \ker\left(H_{n}\left(\vcenter{\xymatrix{B \ar@{=}[r] \ar[d]_-{f} \ar@{}[dr]|{\Rightarrow} & B \ar[d]^-{!_{B}} \\
A \ar[r]_-{!_{A}} & 0}},I_{k+1}\right)\right)\colon K[H_{n}(f,I_{k+1})]\to H_{n}(B,I_{k}).
\]
\end{lemma}
\begin{proof}
This follows from the previous lemma and the naturality of $\gamma^{n}$. Indeed, its naturality square at the map $(1_{B},!_{A})$ is nothing but
\[
\xymatrix{K[H_{n}(f,I_{k+1})] \ar[rrr]^-{\ker H_{n}((1_{B},!_{A}),I_{k+1})} \ar[d]_-{\gamma^{n}_{f}} &&& K[H_{n}(!_{B},I_{k+1})] \ar[d]^-{\gamma^{n}_{!_{B}}}\\
H_{n}(B,I_{k}) \ar@{=}[rrr] &&& H_{n}(B,I_{k});}
\]
and all kernels may be chosen in such a way that $\gamma_{!_{B}}^{n}$ is an identity.
\end{proof}

\begin{proposition}\label{Proposition-Satellites}
Let $I\colon{\A\to\B}$ be a reflector of a semi-abelian category $\A$ with enough projectives onto a Birkhoff subcategory $\B$ of $\A$. Let $n\geq 1$ and $k\geq 0$. Then $H_{n+1}(-,I_{k})\colon{\Extk\A\to\Extk\A}$ with the connecting natural transformation 
\begin{equation}\label{Diagram-Proposition-Satellites}
\vcenter{\xymatrix@!0@=3,5em{& \Extkk\A \ar[ld]_-{\cod} \ar[rd]^-{H_{n}(-,I_{k+1})}\\
\Extk\A \ar@{:>}[rr]^-{\delta^{n+1}} \ar@{.>}[rd]_-{H_{n+1}(-,I_{k})} && \Extkk \A \ar[ld]^-{\ker}\\
& \Extk\A}}
\end{equation}
is the pointwise left satellite of $H_{n}(-,I_{k+1})$. In particular, for any object $A$ of $\A$,
\[
H_{n+1}(A,I) = \Ran_{\cod}(\ker\comp H_n(-,I_{1}))(A) = \lim_{(f,g)\in|(A\downarrow\cod)|}K[H_n(f,I_{1})].
\]
\end{proposition}
\begin{proof}
Let $A$ be an object of $\Extk\A$. Let $p\colon{P\to A}$ be a projective presentation of $A$. We have to show that $(H_{n+1}(A,I_k),\overline{\delta^{n+1}})$ is the limit of 
\[
\xymatrix{(A\downarrow\cod) \ar[rr]^U && \Extkk\A \ar[rr]^{H_n(-,I_{k+1})} && \Extkk\A \ar[rr]^{\ker} &&\Extk\A.}
\]
To do so, let $(L,\lambda)$ be another cone on $\ker\comp H_n(-,I_{k+1})\comp U$; we use the presentation $p$ of $A$ to construct a map of cones $l\colon {L\to H_{n+1}(A,I_{k})}$.

First we consider the case $n=1$. Recall from \cite{Tomasthesis} that by definition $H_1(-,I_m)=I_m$ for all $m\in\N$. Since $p\colon{P\to A}$ is a projective presentation of $A$, and thus $H_2(P,I_{k})=0$, the lower end of the Everaert Sequence~\eqref{Everaert-Sequence} of $p$ becomes
\[
\xymatrix{0 \ar[r] & H_{2}(A,I_{k}) \ar@{{ |>}->}[r]^-{\delta^{2}_{p}} & K[I_{k+1}p] \ar[r]^-{\gamma^{1}_p} & I_{k}P.}
\]
In other words, $\delta^{2}_{p}$ is the kernel of $\gamma^{1}_p$. Recalling Diagram~\eqref{Subdiagram-Exact-Sequences}, consider the following two morphisms in $(A\downarrow \cod)$:
\begin{equation}\label{Diagram-Init}
\vcenter{\xymatrix{P \ar@{ >>}[r]^-p \ar@{=}[d] \ar@{}[dr]|{\Downarrow} & A \ar[d]^{!_A} \ar@{=}[r]^-{1_A} \ar@{}[dr]|{\Downarrow} & A\\
P\ar@{ >>}[r]^{!_P} \ar@{}[dr]|{\Uparrow} & 0 \ar@{}[dr]|{\Uparrow} & A \ar@{=}[u] \ar[l]_-{!_A} \ar@{=}[d] \\
0 \ar[u]^{\gnab_P} \ar@{=}[r]_-{1_0} & 0 \ar[u]_{1_0} & A \ar[l]^-{!_A}}}
\end{equation}
 By Lemma~\ref{Lemma-Gamma}, the naturality of $\lambda$ at the downward-pointing morphism in Diagram~\eqref{Diagram-Init} means $\gamma_p^{1}\comp \lambda_{(p,1_A)}=\lambda_{(!_P,!_A)}$.
This latter morphism is zero, since the naturality of $\lambda$ at the upward-pointing morphism in~\eqref{Diagram-Init} means $\lambda_{(!_P,!_A)}=I_k(\gnab_P)\comp \lambda_{(1_0,!_A)}$, and $I_k0=0$. Hence there exists a unique morphism $l\colon {L\to H_2(A,I_k)}$ satisfying $\lambda_{(p,1_A)}=\delta_p^{2}\comp l$.
 
Higher up in the Everaert Sequence~\eqref{Everaert-Sequence} of $p$, for $n\geq 2$, Corollary~\ref{Corollary-Homology-of-extension} gives us the isomorphism 
\[
\delta^{n+1}_{p}\colon H_{n+1}(A,I_{k})\overset{\cong}{\to} K[H_{n}(p,I_{k+1})].
\]
Here we may simply put $l=(\delta_{p}^{n+1})^{-1}\comp \lambda_{(p,1_A)}$. 

It remains to be shown that, in both cases, the constructed map $l$ is a map of cones. Given any object $(f\colon {B \to C}, g\colon{A\to C})$ of $(A\downarrow \cod)$, there is a map 
\[
\xymatrix{P \ar@{ >>}[r]^-p \ar[d] \ar@{}[dr]|{\Downarrow} & A \ar[d]^-g \ar@{}[dr]|{\Downarrow} & A \ar@{=}[l] \ar@{=}[d]\\
B \ar@{ >>}[r]_-f & C & A \ar[l]^-g}
\]
as $P$ is projective. Writing $h$ for the image of this map under $\ker\comp H_n(-,I_{k+1})\comp U$, we see that the diagram 
\[
\xymatrix@!0@=4,5em{L \ar[rr]^-l \ar[d]_{\lambda_{(p,1_{A})}} \ar[drr]|(.3){\lambda_{(f,g)}}& & H_{n+1}(A,I_k) \ar[dll]|(.3){\delta_p}|(.5){\hole} \ar[d]^{\overline{\delta}_{(f,g)}}\\
K[H_n(p,I_{k+1})] \ar[rr]_{h} &&K[H_n(f,I_{k+1})]}
\]
commutes: $\lambda_{(f,g)}=h\comp \lambda_{(p,1_A)}=h\comp \delta_p \comp l = \overline{\delta}_{(f,g)}\comp l$. Thus $l$ is indeed a map of cones, and $H_{n+1}(A,I_k)$ is the limit of the given diagram.
\end{proof}

\begin{remark}\label{Remark-Symmetry-2-Rest}
This gives a way to derive the $H_{n+1}(-,I_{k})$ from $H_{n}(-,I_{k+1})$ for $n\geq 2$ in exactly the same way as $H_{2}(-,I_{k})$ is derived from $H_{1}(-,I_{k+1})=I_{k+1}$. In other approaches such as \cite{Tomasthesis,EGVdL} the two cases are formally different.
\end{remark}

\subsection{$H_{n+1}(-,I_{k})$ as a satellite of $I_{k+n}$}\label{Subsection-Main-Theorem}

Proposition~\ref{Proposition-Satellites} gives a way to construct $H_{n+1}(-,I_k)$ out of $H_{n}(-,I_{k+1})$. Here, with Theorem~\ref{Theorem-Satellites}, we obtain a one-step construction of $H_{n+1}(-,I_k)$ out of $I_{n+k}$. To be able to apply Proposition~\ref{Proposition-Satellites} repeatedly, we have to show that satellite diagrams like Diagram~\eqref{Diagram-Proposition-Satellites} may be composed in a suitable way (cf.~\cite[Theorem~9]{Janelidze-Satellites}).

The kernel functor
\[
\ker \colon{\Extkk\A \to \Extk\A}
\]
that maps an extension $f\colon {B\to A}$ to its kernel $K[f]$ has a left adjoint, namely the functor ${\Extk\A \to \Extkk \A}$ that sends an object $C$ of $\Extk\A$ to the extension $!_C\colon{C\to 0}$. This allows us to use the following result.

\begin{proposition}\label{Composing-Kan-Extensions}
Suppose that $(I',\delta')=\Ran_{F'}G'I''$ and $(H,\delta)=\Ran_{F}GI'$ as in the diagrams
\[
\vcenter{
\xymatrix@!0@=3,5em{& \A' \ar[ld]_-{F} \ar[rd]^-{I'}\\
\A \ar@{.>}[rd]_-{H} \ar@{:>}[rr]^-{\delta} && \B' \ar[ld]^-{G}\\
& \B }}
\qquad \text{and}\qquad
\vcenter{\xymatrix@!0@=3,5em{& \A'' \ar[ld]_-{F'} \ar[rd]^-{I''}\\
\A' \ar@{.>}[rd]_-{I'} \ar@{:>}[rr]^-{\delta'} && \B''. \ar[ld]^-{G'}\\
& \B'}}
\]
If $G$ is a right adjoint then $(H,G\delta'\comp\delta_{F'})=\Ran_{FF'}GG'I''$: the two diagrams may be composed to form a single Kan extension diagram
\[
\xymatrix@!0@=3,5em{& \A'' \ar[ld]_-{FF'} \ar[rd]^-{I''}\\
\A \ar@{.>}[rd]_-{H} \ar@{:>}[rr]^-{G\delta'\comp\delta_{F'}} && \B''. \ar[ld]^-{GG'}\\
& \B }
\]
If $G$ preserves limits and $(I',\delta')$ and $(H,\delta)$ are pointwise satellites then $(H,G\delta'\comp\delta_{F'})$ is also a pointwise satellite.
\end{proposition}
\begin{proof}
We prove the pointwise case. Let $A$ be an object of $\A$, and $(C,\sigma)$ a cone on the diagram ${GG'I''U\colon{(A\downarrow FF')\to \B}}$. 

For any $A'$ in $\A'$, the pair $(I'A',\overline{\delta'})$ is the limit of the diagram
\[
G'I''U'\colon {(A'\downarrow F')\to \B'}.
\]
Since $G$ preserves limits, $(GI'A',G\overline{\delta'})$ is the limit of $GG'I''U'\colon{(A'\downarrow F')\to \B}$. Now for every $\alpha\colon{A\to FA'}$ the collection $(\sigma_{(A'',F\alpha'\circ \alpha)})_{(A'',\alpha')\in|(A'\downarrow F')|}$ also forms a cone on $GG'I''U'$; hence there is a unique map $\mu_{(A',\alpha)}\colon{C\to GI'A'}$ such that $G\delta'_{A''}\comp GI'\alpha'\comp \mu_{(A',\alpha)}= \sigma_{(A'',F\alpha'\circ \alpha)}$. 

The collection $(\mu_{(A',\alpha)})_{(A',\alpha)\in|(A\downarrow F)|}$ in turn forms a cone on the diagram 
\[
GI'U\colon{(A\downarrow FF')\to \B}.
\]
Indeed, if $(B',\beta)$ is an object of $(A\downarrow F)$ and $f'\colon{B'\to A'}$ is a map in $\A'$ such that $Ff'\comp \beta=\alpha$, then $GI'f'\comp \mu_{(B',\beta)}=\mu_{(A',\alpha)}$, because for every $(A'',\alpha')\in|(A'\downarrow F')|$,
\begin{align*}
G\delta'_{A''}\comp GI'\alpha'\comp GI'f'\comp \mu_{(B',\beta)} &= \sigma_{(A'',F(\alpha'\circ f')\circ \beta)}\\
&=\sigma_{(A'',F\alpha'\circ \alpha)}\\
&=G\delta'_{A''}\comp GI'\alpha'\comp \mu_{(A',\alpha)},
\end{align*}
and the $G\delta'_{A''}\comp GI'\alpha'$ are jointly monic.

This cone gives rise to the needed unique map $c\colon{C\to HA}$. Since it satisfies $\mu_{(A',\alpha)}=\delta_{A'}\comp H\alpha\comp c$ for all $(A',\alpha)\in|(A\downarrow F)|$, we have that
\begin{align*}
G\delta'_{A''}\comp \delta_{F'A''}\comp H\alpha''\comp c 
&= G\delta'_{A''}\comp \delta_{F'A''}\comp HF\alpha'\comp H\alpha\comp c \\
&= G\delta'_{A''}\comp GI'\alpha' \comp \delta_{A'}\comp H\alpha\comp c \\
&= G\delta'_{A''}\comp GI'\alpha' \comp \mu_{(A',\alpha)}\\
&= \sigma_{(A'',F\alpha'\circ \alpha)}=\sigma_{(A'',\alpha'')}
\end{align*}
for all $\alpha''=F\alpha'\comp \alpha\colon{A\to FA'\to FF'A''}$ in $|(A\downarrow FF')|$---and any $\alpha''$ allows such a decomposition.
\end{proof}

\begin{theorem}\label{Theorem-Satellites}
Let $I\colon{\A\to\B}$ be a reflector of a semi-abelian category with enough projectives $\A$ onto a Birkhoff subcategory $\B$ of $\A$. Let $k\geq 0$ and $n\geq 1$. Then 
\[
 H_{n+1}(-,I_k)\colon{\Extk\A\to\Extk\A}
\]
with the connecting natural transformation
\[
\del^{n+1}=\ker^{n-1}\delta^{2}\comp\cdots\comp \ker\delta^{n}\comp\delta^{n+1}\colon{H_{n+1}(-,I_k)\comp \cod^{n}\To \ker^{n} \comp I_{k+n}}
\]
is the pointwise left satellite of $I_{k+n}$.
\[
\xymatrix@!0@=3,5em{& \Extkn\A \ar[ld]_-{\cod^{n}} \ar[rd]^-{I_{k+n}}\\
\Extk\A \ar@{.>}[rd]_-{H_{n+1}(-,I_{k})} \ar@{:>}[rr]^-{\del^{n+1}} && \Extkn\A \ar[ld]^-{\ker^{n}}\\
& \Extk\A}
\]
In particular, for any object $A$ of $\A$,
\[
H_{n+1}(A,I)=\Ran_{\cod^{n}}(\ker^n\comp I_{n})(A)=\lim_{(f,g)\in|(A\downarrow \cod^{n})|}K^n[I_{n}f].
\]
\end{theorem}
\begin{proof}
This follows from gluing diagrams as in Proposition~\ref{Proposition-Satellites} together using Proposition~\ref{Composing-Kan-Extensions}.
\end{proof}

\subsection{Symmetry}\label{Subsection-Gamma}

Proposition~\ref{Proposition-Satellites} gives an interpretation of the connecting morphisms $\delta^{n}_{f}$ in the Everaert sequence as left satellites. The connecting morphisms $\gamma^{n}_{f}$ have a dual interpretation: $(H_{n}(-,I_{k}),\gamma^{n})$ is a right satellite (left Kan extension) of $H_{n}(-,I_{k+1})$. 

\begin{proposition}\label{Proposition-Satellites-Ext-Dual}
Let $I\colon{\A\to\B}$ be a reflector of a semi-abelian category $\A$ with enough projectives onto a Birkhoff subcategory $\B$ of $\A$. Consider $n\geq 1$ and $k\geq 0$. Then $(H_{n}(-,I_{k}),\gamma^{n})$, i.e., $H_{n}(-,I_{k})\colon{\Extk\A\to\Extk\A}$ with the connecting natural transformation 
\[
\xymatrix@!0@=3,5em{& \Extkk\A \ar[ld]_{H_{n}(-,I_{k+1})} \ar[rd]^{\dom}\\
\Extkk\A \ar[rd]_{\ker} \ar@{:>}[rr]^-{\gamma^{n}} && \Extk\A \ar@{.>}[ld]^{H_{n}(-,I_{k})}\\
& \Extk\A}
\]
is the pointwise right satellite of $H_{n}(-,I_{k+1})$.
\end{proposition}
\begin{proof}
For any $A$, the category $(\dom\!\downarrow\!A)$ has a terminal object $(!_{A}\colon A\to 0,1_A)$, so the colimit object of the diagram 
\[
\xymatrix@!@=3,5em{(\dom\!\downarrow\!A)\ar[r]^-U & \Extkk \A \ar[r]^-{H_{n}(-,I_{k+1})} & \Extkk\A \ar[r]^-{\ker} & \Extk\A}
\]
is $K[H_{n}(!_A,I_{k+1})]=H_{n}(A,I_{k})$. The component of the colimit cocone at
\[
(g\colon{B\to C},f\colon{B\to A})\in|(\dom\!\downarrow\!A)|
\]
is 
\[
\resizebox{\textwidth}{!}{\mbox{$
\begin{aligned}
\ker\left(H_{n}\left(\vcenter{\xymatrix{B \ar[r]^{f} \ar[d]_-{g} \ar@{}[dr]|{\Rightarrow} & A \ar[d]^-{!_{A}} \\
C \ar[r]_-{!_{C}} & 0}},I_{k+1}\right)\right)
&=\ker\left(H_{n}\left(\vcenter{\xymatrix{B \ar[r]^-{f} \ar[d]_-{!_{B}} \ar@{}[dr]|{\Rightarrow} & A \ar[d]^-{!_{A}} \\
0 \ar@{=}[r] & 0}},I_{k+1}\right)\right)
\comp 
\ker\left(H_{n}\left(\vcenter{\xymatrix{B \ar@{=}[r] \ar[d]_-{g} \ar@{}[dr]|{\Rightarrow} & B \ar[d]^-{!_{B}} \\
C \ar[r]_-{!_{C}} & 0}},I_{k+1}\right)\right)\\
&=H_{n}(f,I_{k})\comp \gamma^{n}_{g}\\
&=\overline{\gamma^{n}}_{(g,f)}
\end{aligned}
$}}
\]
by Lemma~\ref{Lemma-Gamma} and Lemma~\ref{Lemma-iota}.
\end{proof}

\section{Homology without projectives}\label{Section-Without-Projectives}

In this section we set up a homology theory without projectives by defining homology via pointwise satellites as they appear in Proposition~\ref{Proposition-Satellites}. 

\begin{proposition}\label{Proposition-Equivalence-Homology-Definitions}
Let $\A$ be a pointed exact protomodular category and $I\colon{\A\to\B}$ a reflector onto a Birkhoff subcategory $\B$ of $\A$. Consider $k\geq 0$ and $A\in|\Extk\A|$. If it exists, write
\[
H_{(2,k)}=\Ran_{\cod} (\ker\comp I_{k+1}) 
\]
for the pointwise left satellite of $I_{k+1}$ relative to the functors $\cod$ and $\ker$. Now suppose $H_{(n,k+1)}$ exists for $n\geq 2$, and write
\[
H_{(n+1,k)}=\Ran_{\cod} (\ker\comp H_{(n,k+1)})                                                                              \]
for the pointwise left satellite of $H_{(n,k+1)}$ relative to $\cod$ and $\ker$, if this exists. Then $H_{(n+1,k)}$ is also the left satellite of $I_{k+n}$ relative to the functors $\cod^n$ and $\ker^n$.
\end{proposition}
\begin{proof} 
The proof is the same as the proof of Theorem~\ref{Theorem-Satellites}. 
\end{proof}

\begin{definition}\label{Definition-Homology}
Let $\A$ be a pointed exact protomodular category and $I\colon{\A\to\B}$ a reflector onto a Birkhoff subcategory $\B$ of $\A$. Consider $k\geq 0$ and $A\in|\Extk\A|$, and let $n\geq 1$. If the functor $H_{(n+1,k)}$ from Proposition~\ref{Proposition-Equivalence-Homology-Definitions} exists, we call it the \defn{$(n+1)$-st homology functor
\[
H_{n+1}(-,I_{k})\colon{\Extk\A\to\Extk\A}
\]
(with coefficients in $I_{k}$)}.
\[
\xymatrix@!0@=3,5em{& \Extkk\A \ar[ld]_-{\cod} \ar[rd]^-{H_{n}(-,I_{k+1})}\\
\Extk\A \ar@{:>}[rr]^-{\delta^{n+1}} \ar@{.>}[rd]_-{H_{n+1}(-,I_{k})} && \Extkk \A \ar[ld]^-{\ker}\\
& \Extk\A}
\]
We also write $H_1(-,I_k)=I_k$.
\end{definition}

\begin{remark}
 From now on, when we write $H_{n+1}(-,I_k)$ we mean the homology functor as defined here via Kan extensions, rather than the homology defined via Hopf formulae as in Sections~\ref{Section-Preliminaries} and~\ref{Section-Satellites}.
\end{remark}

\begin{remark}\label{Remark-Homology}
For any object $A\in|\Extk\A|$, if  $H_{2}(A,I_{k})$ exists, it is the limit object of the diagram
\[
\xymatrix@!@=3,5em{(A\downarrow \cod) \ar[r]^U & \Extkk\A \ar[r]^-{I_{k+1}} & \Extkk\A \ar[r]^-{\ker} &\Extk\A.}
\]
Similarly, if $H_{n+1}(A,I_{k})$ exists, it is the limit object of the diagram
\[
\xymatrix@!@=3,5em{(A\downarrow \cod) \ar[r]^U & \Extkk\A \ar[r]^{H_n(-,I_{k+1})} & \Extkk\A \ar[r]^{\ker} & \Extk\A}
\]
or equivalently of
\begin{equation}\label{Diagram-All-In-One}
%\vcenter
{\xymatrix@1@!@=3,5em{(A\downarrow \cod^n) \ar[r]^U & \Extkn\A \ar[r]^{I_{k+n}} & \Extkn\A \ar[r]^{\ker^n} & \Extk\A.}}
\end{equation}
Potentially, these limits may exist for a given object $A$ even if the homology functors $H_{n+1}(-,I_k)$ do not exist in full. Such a limit is most easily computed pointwise (in $\Arrk\A$) and then shown to be an extension.
\end{remark}

\begin{example}[When the reflection is the identity]\label{Example-Trivial-Birkhoff-Subcategory}
If $\B=\A$ then all $I_{n}$ are identity functors, and the $H_{n}$ are zero for $n\geq 2$. To see this, we have to prove that the functor $0\colon{\Extk\A\to\Extk\A}$ is a pointwise Kan extension of $\ker\colon{\Extkk\A\to \Extk\A}$ along $\cod\colon{\Extkk\A\to \Extk\A}$, for all $k\geq 0$. This shows that $H_{2}$ is zero, which immediately implies that the higher homologies are also zero, being satellites of the zero functor.

Let $A$ be an object of $\Extk\A$ and $({L},\lambda)$ a cone on $\ker\comp U\colon{(A\downarrow \cod)\to \Extk\A}$. Then any map $\lambda_{(f,g)}$, where $(f\colon{B\to C},g\colon{A\to C})\in |(A\downarrow \cod)|$, fits into the commutative diagram
\[
\xymatrix@!0@=4,5em{L \ar[r]^-{\lambda_{(f,g)}} \ar[d]_-{\lambda_{(!_{0},!_{A})}} \ar[rd]|{\lambda_{(!_{B},!_{A})}} & K[f] \ar@{{ |>}->}[d]^-{\Ker f}\\
0 \ar[r] & B,}
\]
which means that $\lambda_{(f,g)}$ is the zero map. If now $(L,\lambda)$ is a limit cone, this implies that $L$ is zero.
\end{example}

The category $(A\downarrow \cod)$ is rather large, and in a given situation it may be very hard to decide whether the needed limits do indeed exist. Even if they do, they may still be hard to compute. But we may replace the above diagrams with simpler ones, for example using the concept of \emph{initial subcategory}. Recall its definition as it occurs in \cite[Section~IX.3]{MacLane}:

\begin{definition}\label{Definition-Initial-Diagram}
An \defn{initial functor} is a functor $F\colon{\D \to \C}$ such that for every object $C$ of $\C$, the comma category $(F\downarrow C)$ is non-empty and connected. A subcategory $\D$ of a category $\C$ is called \defn{initial} when the inclusion of $\D$ into $\C$ is an initial functor, i.e., for every object $C\in|\C|$, the full subcategory $(\D\downarrow C)$ of  $(\C\downarrow C)$ determined by the maps ${D\to C}$ with domain $D$ in $\D$ is non-empty and connected.
\end{definition}

If $\D$ is initial in $\C$ then limits of diagrams over $\C$ may be computed as the limit of their restriction to $\D$. More generally, if $F\colon{\D\to \C}$ is initial then a diagram $G\colon{\C\to \E}$ has a limit if and only if so does $GF$, in which case it may be computed as the limit of $GF$. 

For any object $A$ of $\Extk\A$, let $\ExtkkA\A$ denote the \defn{category of extensions of~$A$}, the pre\-image in $\Extkk\A$ of the arrow $1_A$ under the functor $\cod\colon {\Extkk\A\to \Extk\A}$. Then the functor $U'\colon{\ExtkkA\A\to (A\downarrow \cod)}$ that sends an extension $f\colon{B\to A}$ of $A$ to the pair $(f, 1_{A})$ is easily seen to be initial: for every object
\[
(f\colon{B\to C}, g\colon{A\to C})
\]
of $(A\downarrow \cod)$ there is the natural morphism ${U'\overline f \to (f,g)}$
\[
\xymatrix{\overline{B} \ar[d] \ar@{ >>}[r]^-{\overline{f}}\ar@{}[rd]|<{\pullback} \ar@{}[dr]|{\Downarrow} & A \ar[d]^-g \ar@{=}[r]^-{1_A} \ar@{}[dr]|{\Downarrow} & A \ar@{=}[d] \\
B \ar@{ >>}[r]_-f & C & A, \ar[l]^-g}
\]
where $\overline f$ is the pullback of $f$ along $g$; this $\overline{f}$ is an extension by \cite[Proposition~3.5]{EGVdL}. Also, any other morphism 
\[
\xymatrix{D \ar[d] \ar@{ >>}[r]^-{h} \ar@{}[dr]|{\Downarrow} & A \ar[d]^-g \ar@{=}[r]^-{1_A} \ar@{}[dr]|{\Downarrow} & A \ar@{=}[d] \\
B \ar@{ >>}[r]_-f & C & A, \ar[l]^-g}
\]
factors over this morphism ${U'\overline f \to (f,g)}$, by the universal property of a pullback.
This means that the limit of $\ker\comp H_n(-,I_{k+1})\comp U$ may also be computed as the limit of $\ker\comp H_n(-,I_{k+1})\comp UU'$ and moreover, since $U U'$ is just the inclusion of the subcategory $\ExtkkA\A$ into $\Extkk\A$, as the limit of
\[
\ker\comp H_n(-,I_{k+1})\colon{\ExtkkA\A\to \Extk\A.}
\]

But even now the diagram of shape $\ExtkkA \A$ over which the limit is computed may be too large, in the sense that even if $\A$ is small-complete, it is still unclear whether the limit of $\ker\comp H_n(-,I_{k+1})$ exists. In the case where $\A$ has enough projectives, however, it is possible to further cut down on the size of this diagram. In this case Proposition~\ref{Proposition-Satellites} shows that the limit of this diagram exists and is equal to the homology object defined via the Hopf formulae. But making the diagram smaller gives a new way to calculate this homology. This situation is discussed in Section~\ref{Section-Homology-With-Projectives}.

\begin{notation}\label{Notation-Extensions}
Let $A\in|\Extk\A|$. Denote by $\ExtknA\A$ the category of \defn{$n$-extensions of $A$}, defined as the preimage of the arrow $1_A$ under the functor
\[
\cod^n\colon {\Extkn\A\to \Extk\A}.
\]
This generalises the category $\ExtkkA\A$ of \emph{extensions of $A$} defined above. Thus the objects are $n$-extensions with ``terminal object'' $A$, when viewed as diagrams in the category $\Ext^k\!\A$, and the maps are those maps in $\Extkn\A$ which restrict to the identity on $A$ under $\cod^n$. Similarly the category $\CExt^{k+n}_{A}\!\A$ denotes the full subcategory of $\ExtknA\A$ determined by those $n$-extensions which are central. The Birkhoff subcategory $\B$ is understood, and not mentioned in the notation.
\end{notation}

\begin{remark}\label{Remark-Notation-Extensions}
The functor $U'\colon {\ExtknA\A\to (A\downarrow \cod^{n})}$ which sends an $n$-extension $f$ of $A$ to $(f,1_A)$ is still initial. This may be shown by induction, using the fact that in a category of $n$-fold extensions, the $(n+1)$-extensions are pullback-stable \cite[Proposition~3.5]{EGVdL}.
\end{remark}

\begin{proposition}\label{Proposition-One-Step-Homology}
Let $I\colon{\A\to\B}$ be a reflector of a pointed exact protomodular category $\A$ onto a Birkhoff subcategory $\B$ of $\A$. Consider $k\geq 0$ and $n\geq 1$ and $A\in |\Extk\A|$. If it exists, $H_{n+1}(A,I_k)$ is also the limit of the diagram
\[
\ker^n\comp I_{k+n}\colon{\Ext^{k+n}_A\!\A\to \Extk\A.}
\]
\end{proposition}
\begin{proof} 
This uses Diagram~\eqref{Diagram-All-In-One} and the fact that $U'\colon{\ExtknA\A \to (A\downarrow \cod^n)}$ is initial.
\end{proof}

\begin{corollary}\label{Corollary-Diagram-Central-Extensions}
For $k\geq 0$, $n\geq 1$ and $A\in|\Extk\A|$, if it exists, $H_{n+1}(A,I_{k})$ is the limit of the diagram
\[
\ker^{n}\colon{\CExt^{k+n}_A\!\A\to \Extk\A.}
\]
\end{corollary}
\begin{proof}
The functor $I_{k+n}\colon{\Ext^{k+n}_A\!\A\to \CExt^{k+n}_A\!\A}$ is initial because, for any central extension ${f\in|\CExt^{k+n}_A\!\A|}$, we have $I_{k+n}f=f$, so the comma category $(I_{k+n}\downarrow f)$ is non-empty and connected.
\end{proof}

Since limits commute with kernels, Corollary~\ref{Corollary-Diagram-Central-Extensions} also says that $H_{n+1}(A,I_{k})$ may be computed as the $n$-fold kernel of a certain $(n+k)$-fold arrow in $\A$, namely, the limit in $\Arr^{k+n}\!\A$ of the inclusion of $\CExt^{k+n}_A\!\A$ into $\Arr^{k+n}\!\A$. Sometimes this $n$-fold arrow in $\Arrk\A$ itself happens to be an $n$-fold central extension of $A$. We say that an $n$-fold central extension of a $k$-extension $A$ is \defn{universal} when it is an initial object of $\CExt^{k+n}_{A}\!\A$. Recall from~\cite{Gran-VdL} (but see also \cite{Froehlich, Janelidze-Kelly}) that, when $\A$ is a semi-abelian category and $I=\ab\colon{\A\to \Ab\A}$ is the abelianisation functor, then an object $A$ of $\A$ admits a universal central extension $p$ if and only if it is \defn{perfect}: its abelianisation is zero. In this case, $H_2(A,\ab)$ is the kernel of $p$. This latter property holds in general, also for higher extensions:

\begin{corollary}\label{Corollary-Initial-Central-Extension}
Consider $k\geq 0$, $n\geq 1$ and $A\in|\Extk\A|$. If $A$ has a universal $n$-fold central extension $p$ then $H_{n+1}(A,I)=K^n[p]$. In particular, if $A\in|\A|$ has a universal central extension $p\colon{P\to A}$ then $H_2(A,I)=K[p]$.
\end{corollary}
\begin{proof}
The limit of a functor over a category that has an initial object is the value of the functor at this object.
\end{proof}

\begin{example}[The homology of zero is zero]\label{Example-Zero}
If $A=0$ then, for any $n\geq 1$, the category $\CExtnA\A$ has an initial object, the zero $n$-cube. Taking kernels as in Corollary~\ref{Corollary-Initial-Central-Extension} gives $H_{n+1}(0,I)=0$.
\end{example}

\begin{remark}\label{Remark-Identity-Weakly-Universal}
Note that in certain special cases a \emph{weakly} universal extension can also determine the homology of a $k$-extension $A$. When $1_A$ is a weakly universal extension of $A$, i.e., if every extension $f\colon {B\to A}$ of $A$ is split, we have $H_2(A,I_k)=0$. This is because $K[I_k1_A]=0$ for any object $A$, so if $1_A$ is weakly initial, every leg of a cone over $\ker\comp I_k\colon {\ExtkkA\A\to \Extk\A}$ factors over $K[I_k1_A]$ and thus is zero. In particular, we get:
\end{remark}

\begin{example}[The homology of a projective object is zero]
For any projective object $P$ and any $n\geq1$ we have $H_{n+1}(P,I)=0$, since $1_P$ (and also the $n$-extension only consisting of the maps $1_P$) is always weakly initial when $P$ is projective.
\end{example}

\begin{example}[Homology of finite groups]\label{Example-Groups}
For a finite group, we compare its second homology groups with respect to two different adjunctions. On the one hand we have the abelianisation functor $\ab\colon {\Gp \to \AbGp}$, where $\Gp$ is the category of groups, $\Ab\Gp$ is the Birkhoff subcategory of abelian groups, and $\ab G=G/[G,G]$. This example has been studied in the classical setting in~\cite{EverVdL2} (for lower dimensions) and in~\cite{Tomasthesis, EGVdL} (higher dimensions). Here the centralisation functor $\ab_1$ takes an extension $f\colon {B \to A}$ to $\centr f\colon {B/[K[f],B] \to A}$. As mentioned in Section~\ref{Subsection-Barr-Beck-Definition}
, in this case Definition~\ref{Definition-Homology} 
gives the classical integral homology of groups. 

On the other hand, we could focus on finite groups and let $\A=\Fin\Gp$ be the category of finite groups and $\B=\Fin\Ab=\Ab\Fin\Gp$ its Birkhoff subcategory of finite abelian groups. Note that $\Fin\Gp$ is not semi-abelian and doesn't have enough projectives, but nevertheless it is pointed, Barr exact and Bourn protomodular. Here $I\colon {\A \to \B}$ again sends a group $G$ to $\finab G=G/[G,G]$ and $I_1\colon {\Ext\Fin\Gp \to \Ext\Fin\Gp}$ sends an extension $f\colon {B\to A}$ to
\[
\fincentr f\colon {\frac{B}{[K[f],B]}\to A}.
\]
We show that, for any finite group, its second homology groups with respect to the two theories coincide.

For perfect groups this is clear. Recall from Corollary~\ref{Corollary-Initial-Central-Extension} 
that if a group $G$ has a universal central extension $p\colon {P\to G}$, then the homology is $H_2(G,\ab)=K[p]$; this is the case when $G$ is perfect: $\ab G=0$. So given a finite perfect group $G$, we know that it has a universal central extension $p\colon {P\to G}$ in the category $\Gp$ of all groups, and that $H_2(G,\Z)=H_2(G,\ab)=K[p]$. But we also know that the integral homology of a finite group is a \emph{finite} group, therefore the group $P$ must also be finite, and the universal central extension $p\colon {P\to G}$ lies in the category $\Fin\Gp$ of \emph{finite} groups. Thus we also have $H_2(G,\finab)=K[p]$. So for a finite perfect group $G$ we have $H_2(G,\finab)=H_2(G,\ab)=H_2(G,\Z)$.

For a general group, we need a few more steps to prove this equality. 
\vspace{\baselineskip}

\paragraph*{\textit{Step 1:}}

First we want to show that, for any finite group $G$, there is a central extension $G^*\to G$ with kernel $H_2(G,\Z)$, such that in the diagram 
\begin{equation}\label{big-diagram}
\ker\colon {\CExt_{G}\Gp \to \Gp},
\end{equation}
the leg from the limit $H_2(G,\ab)$ to this object is an isomorphism. We consider \defn{stem extensions}: central extensions $g\colon {H\to G}$ with $K[g]\leq [H,H]$. This condition implies that $\ab H \to \ab G$ is an isomorphism, or equivalently that the map $K[g]\to \ab H$ is zero. So it follows from exactness in \eqref{Everaert-Sequence} that the leg $H_2(G,\ab)\to K[g]$ is a surjection when $g$ is a stem extension. To find a stem extension with $H_2(G,\Z)$ as its kernel, we use the \defn{Schur multiplier} $M(G)$ of a finite group $G$ introduced in \cite{Schur1904}. Schur proved in \cite{Schur1907} that for a finite group $G$, this  multiplier $M(G)$ may be expressed in terms of what is now called the Hopf formula (which, in the infinite case, was only introduced in~\cite{Hopf}), and so we have $M(G)\iso H_2(G,\Z)$ (see also, e.g., \cite[Theorem 2.4.6]{Karp:Schur-Multiplier}). In \cite{Schur1904} he showed that, for any finite group $G$, there is a stem extension $f\colon {G^*\to G}$ of $G$ with kernel $M(G)$ (see also \cite[Theorem 2.1.4]{Karp:Schur-Multiplier}). 

Putting these two facts together, we see that $H_2(G,\Z)$ occurs in the diagram \eqref{big-diagram} as the kernel of this stem extension $f$, and that the leg from $H_2(G,\ab)$ to it must be an isomorphism, being a surjection between finite groups of the same size. From now on we shall assume that this isomorphism is an identity.
\vspace{\baselineskip}

\paragraph*{\textit{Step 2:}}

We now consider the diagram of kernels of \emph{finite} central extensions of $G$,
\begin{equation}\label{finite-diagram}
\ker\colon {\CExt_{G}\Fin\Gp \to \Fin\Gp,}
\end{equation}
which is a small diagram and so has a limit in $\Gp$ which we denote by $L$. We shall show in Step 3 that $L\iso H_2(G,\ab)$ and so is actually the limit of \eqref{finite-diagram} in the category $\Fin\Gp$ as well, as $H_2(G,\ab)$ is a finite group.

$H_2(G,\ab)$ forms a cone on \eqref{finite-diagram}, using the legs from \eqref{big-diagram}. The induced map of cones to $L$ gives a splitting for the leg $p\colon {L\to H_2(G,\Z)=K[f]}$. As these are all abelian groups, we have $L\iso H_2(G,\Z)\oplus E$ for some abelian group $E$, and $p=\pi_1\colon {L\to H_2(G,\Z)}$, the first projection. We consider the following central extensions and maps between them:
\[
\xymatrix{H_2(G,\Z) \ar@{ |>->}[r] \ar@<.5ex>@{ |>->}[d] & G^* \ar@{ >>}[r]^-f \ar@<.5ex>@{ |>->}[d]^-{(1_{G^*},0)} & G \ar@{=}[d]\\
H_2(G,\Z)\oplus E \ar@<.5ex>@{ >>}[u]^p \ar@{ |>->}[r] \ar@{ >>}[d]_-{\pi_2} & G^*\times E \ar@<.5ex>@{ >>}[u]^{\pi_1} \ar@{ >>}[r]^-{f\circ\pi_1} \ar@{ >>}[d]^-{f\times 1_E} & G \ar@{=}[d]\\
E \ar@{ |>->}[r] & G\times E \ar@{ >>}[r]_-{\pi_1} & G}
\]
Since the extension $\pi_1 \colon {G\times E \to G}$ is split, the leg from $L$ to $E=K[\pi_1]$ must be the zero map. So the leg from $L$ to $K[f\comp\pi_1]$ is
\[
1_{H_2(G,\Z)}\oplus 0 \colon {L\iso H_2(G,\Z)\oplus E \to H_2(G,\Z)\oplus E},
\]
as $H_2(G,\Z)\oplus E$ is a product. 
\vspace{\baselineskip}

\paragraph*{\textit{Step 3:}}

Finally we consider a third, even smaller diagram. Let $\C$ be the full subcategory of $\CExt\Fin\Gp$ containing those extensions $g$ of $G$ for which there exists a map $f\to g$ in $\CExt_{G}\Fin\Gp$. We consider the subdiagram
\begin{equation}\label{restricted-diagram}
\ker\colon {\C\to \Fin\Gp,}
\end{equation}
the limit of which is $H_2(G,\ab)$. For any cone $D$ over the diagram \eqref{restricted-diagram}, the two legs $d\colon {D\to H_2(G,\Z)=K[f]}$ and $0\colon {D\to E=K[\pi_1]}$ again determine the leg to the product, $(d,0)\colon {D\to H_2(G,\Z)\oplus E=K[f\comp\pi_1]}$. The leg $d$ also forms the unique cone map ${D\to H_2(G,\ab)}$. Notice that in \eqref{restricted-diagram} we also have maps from $H_2(G,\Z)\oplus E$ to any other object, as $p\colon {H_2(G,\Z)\oplus E \to H_2(G,\Z)}$ is part of the diagram. So as we have $(1_{H_2(G,\Z)}\oplus 0)\comp (d,0)=(d,0)$, the map $(d,0)\colon {D\to L}$ is a cone map and makes $L$ into a limit of \eqref{restricted-diagram}. So $L\iso H_2(G,\ab)$ as promised, and we have $H_2(G,\finab)=H_2(G,\ab)=H_2(G,\Z)$ for any finite group $G$. 
\end{example}

\begin{example}[Internal groups in an exact category]\label{Example-Internal-Groups-Exact-Category}
A possible source of further examples is the category of internal groups $\Gp\E$ in an exact category $\E$, with its Birkhoff subcategory of internal abelian groups $\AbGp\E$. When $\E$ is exact, $\Gp\E$ is semi-abelian if and only if it has coproducts (see~\cite{Janelidze-Marki-Tholen}); it is always pointed exact protomodular. But in general it need not have enough projectives, so our definition of homology via Kan extensions could be a useful tool. One particular class of examples amongst these are internal groups in a topos. The category of internal abelian groups $\AbGp\E$ in a Grothendieck topos $\E$ has enough injectives (see, e.g., Chapter~8 of~\cite{Johnstone:Topos-Theory}), so cohomology theory is possible in this category, but enough projectives are not generally available. In future work we intend to investigate the category of group-valued sheaves on a space as an example of such a situation. Other interesting Birkhoff subcategories of $\Gp\E$ might exist, giving further situations where our definition of homology could be used. 
\end{example}

It is well known that all integral homology groups of a group are abelian. More generally, both approaches to homology discussed in Subsections~\ref{Subsection-Barr-Beck-Definition} and~\ref{Subsection-Everaert-Definition} are such that the homology objects are abelian objects of the Birkhoff subcategory $\B$. We now prove that our homology objects $H_{n+1}(A,I)$ also satisfy these properties. 

\begin{lemma}\label{Lemma-In-B}
Consider an object $A\in|\A|$. The kernel $K[f]$ of a central extension $f\colon {B\to A}$ of $A$ is an object of the Birkhoff subcategory $\B$. More generally, for $k\geq 1$ and $A\in|\Extk\A|$, the kernel $K[f]$ of a central extension $f$ of $A$ is a $k$-fold central extension.
\end{lemma}
\begin{proof}
Let $k\geq 0$, and $A\in|\Extk\A|$. First consider a trivial extension $f\colon {B\to A}$. This means $f$ is the pullback of $I_kf\colon {I_kB\to I_kA}$ along $\eta^k_A$, so $K[f]$ is isomorphic to $K[I_kf]$. This kernel of the extension $I_kf\colon {I_kB\to I_kA}$ is a $k$-fold central extension (or an object of $\B$ for $k=0$) because the category $\CExt^k\A$ is closed under limits which exist in $\Ext^k\A$, as it is a full replete reflective subcategory. (For $k=0$ just note that the Birkhoff subcategory $\B$ is closed under subobjects.)
Now for a central extension $f\colon {B\to A}$, recall from Remark~\ref{Remark-Central-Extensions} that there exists an extension $g$ such that the pullback $\overline f$ of $f$ along $g$ is trivial.
 \[
 \xymatrix{K[\overline f] \ar@{ |>->}[r] \ar@{=}[d] & \overline  B\ar@{ >>}[r]^{\overline f} \ar@{ >>}[d] \ar@{}[dr]|<{\pullback} & \overline A\ar@{ >>}[d]^g \\
K[f] \ar@{ |>->}[r] & B\ar@{ >>}[r]_f &A}
\] 
But then $K[f]=K[\overline f]$, which is a $k$-fold central extension (or an object of $\B$) as $\overline f$ is trivial.
\end{proof}

\begin{remark}
The converse implication does not hold, as for example in the category of groups not every extension with abelian kernel is central.
\end{remark}

\begin{proposition}\label{Proposition-Homology-In-B}
Let $A$ be an object of $\A$ and $n\geq 0$. Then $H_{n+1}(A,I)$ is an object of $\B$.
\end{proposition}
\begin{proof}
If $n=0$ the result is clear as $H_1(A,I)=IA$. For $n\geq 1$, we use Lemma~\ref{Lemma-In-B} repeatedly to see that the diagram from Corollary~\ref{Corollary-Diagram-Central-Extensions} factors over $\B$ and becomes the functor $\ker^{n}\colon{\CExt^{n}_A\A\to \B}$. Since $\B$ is closed under limits in $\A$, the limit $H_{n+1}(A,I)$ of this diagram is still an object of $\B$.
\end{proof}

\begin{example}[When the reflection is zero]\label{Example-Zero-Birkhoff-Subcategory}
If $\B=0$, the zero subcategory in $\A$, then all homology objects are zero, because they are in $\B$ by Proposition~\ref{Proposition-Homology-In-B}.
\end{example} 

The proofs of the next result---Proposition~\ref{Proposition-Homology-Abelian}---and its lemma were offered to us by Tomas Everaert. Recall that an object $A$ of a pointed exact protomodular category $\A$ is \defn{abelian} if it carries an internal abelian group structure. Such a structure is necessarily unique, and is given by a morphism  $m\colon{A\times A\to A}$ satisfying $m\comp (1_A,0) = 1_A = m\comp (0,1_A)$, called its \defn{addition} (see~\cite{Borceux-Bourn}). The abelian objects form a Birkhoff subcategory $\Ab\A$ of $\A$.

\begin{lemma}\label{Lemma-Homology-Abelian}
For any $k\geq 0$ and any $(k+1)$-extension $f\colon{B\to A}$ in $\A$, the image of the connecting morphism 
\[
\delta^2_f\colon{H_2(A,I_k)\to K[H_1(f,I_{k+1})]=K[I_{k+1}f]}
\]
is an abelian object of $\Arrk\A$.
\end{lemma}
\begin{proof}
We show that $\Im[\delta^2_f]$ is a subobject of an abelian object in $\Arrk\A$, namely the kernel of the map
\[
K[(I_{k+1}f, I_k f)]\colon {K[I_{k+1}f]\to K[I_k f]};
\]
here $(I_{k+1}f,I_kf)\colon {\eta_{I_{k+1}[f]} \to \eta_A}$ is a double extension in $\Extk\A$, so its kernel is an extension by~\cite[Proposition 3.9]{EGVdL}. To see that the kernel $K^2[(I_{k+1}f,I_kf)]$ of this extension is an abelian object of $\Arrk\A$, write $(\pi_1,\pi_2)\colon{R[f]\to B\times B}$ for the kernel pair of $f$, and recall the construction of $J_{k+1}[f]$ from Subsection~\ref{Subsection-Reflectors}. We have $K[I_{k+1}f]=K[f]/J_{k+1}[f]=K[f]/\pi_2(J_kR[f]\cap K[f])$, since $J_{k+1}[f]=J_kR[f]\cap K[f]$ as a normal subobject of $R[f]$, and its direct image under $\pi_2$ gives us a normal subobject of $B$ (note that $\pi_2(J_{k+1}[f])=J_{k+1}[f]$ as $\mu^1_f$ is a normal monomorphism). Similarly $K[I_{k}f]=K[f]/(J_kB\cap K[f])=K[f]/(\pi_2J_kR[f]\cap \pi_2K[f])$, so that 
\[
K^2[(I_{k+1}f, I_k f)]=\frac{\pi_2J_kR[f]\cap \pi_2K[f]}{\pi_2(J_kR[f]\cap K[f])}
\]
by Noether's First Isomorphism Theorem~\cite[Theorem~4.3.10]{Borceux-Bourn}. Theorem~2.1 in~\cite{Bourn:Direct-Image} implies that this object is abelian.

Now consider  the arrow $(1_B,!_A)\colon{f\to !_B}$  in $\Extk\A$
\[
\vcenter{\xymatrix{B \ar@{-{ >>}}[d]_-{f} \ar@{=}[r] \ar@{}[dr]|{\Rightarrow} & B \ar@{-{ >>}}[d]^{!_B} \\ 
A \ar[r]_{!_A} & 0}}
\qquad\qquad
\vcenter{\xymatrix{H_2(A,I_k) \ar[d]_-{\delta^2_f} \ar[r] & H_2(0,I_k) \ar[d]^{\delta^2_{!_B}} \\ 
K[I_{k+1}f] \ar[r] & I_kB}}
\]
and the induced commutative square on the right hand side. As $H_2(0,I_k)$ is zero, the map $\delta^2_f$ factors over the kernel of $\eta_{I_{k+1}[f]}\comp\Ker I_{k+1}f\colon{K[I_{k+1}f] \to I_kB}$. The image of this latter map is $K[I_k f]$, so $\Im[\delta^2_f]$ is indeed a subobject of the abelian object $K^2[(I_{k+1}f,I_kf)]$, and thus itself an abelian object as claimed.
\end{proof}

\begin{proposition}\label{Proposition-Homology-Abelian}
Let $A$ be an object of $\A$ and $n\geq 1$. Then $H_{n+1}(A,I)$ is an abelian object of~$\A$.
\end{proposition}
\begin{proof}
It suffices to show that, for all $k\geq 0$ and any $k$-extension $A$, the object $H_{2}(A,I_k)$ is abelian in $\Arrk\A$, as then the higher homology objects are limits of a diagram of abelian objects, and thus abelian by induction. To show $H_2(A,I_k)$ is abelian, consider the functor
\[
H_2(-,I_k)\times H_2(-,I_k)\colon\Extk\A\to\Extk\A 
\]
that sends a $k$-extension $A$ to the product $H_2(A,I_k)\times H_2(A,I_k)$. The previous lemma gives rise to a natural transformation
\[
{(H_2(-,I_k) \times H_2(-,I_k)) \circ \cod  \To \ker \circ I_{k+1}}
\]
of functors from $\Extkk\A$ to $\Extk\A$; the component of this natural transformation at a $(k+1)$-extension $f\colon{B\to A}$ is the composition
\[
H_2(A,I_k) \times H_2(A,I_k) \to \Im[\delta^2_f] \times \Im[\delta^2_f] \to \Im[\delta^2_f] \to K[I_{k+1}f].
\]
Here the first arrow is the corestriction of $\delta^2_f\times \delta^2_f$, the second arrow is the addition on the abelian object $\Im[\delta^2_f]$, and the last arrow is the inclusion of the image into the codomain of $\delta^2_f$. The universal property of the Kan extension $(H_2(-,I_k),\delta^{2})$ now yields a natural transformation ${H_2(-,I_k)\times H_2(-,I_k)\To H_2(-,I_k)}$ which is easily seen to define an abelian group structure on all $H_2(A,I_k)$.
\end{proof}

\section{Homology with projectives}\label{Section-Homology-With-Projectives}

In this section we investigate our new definition of homology in the situation when $\A$ does have enough projectives. In this case we know that homology exists, for example via Everaert's definition using the Hopf formulae, and Proposition~\ref{Proposition-Satellites} shows that it coincides with the notion introduced in Definition~\ref{Definition-Homology}. But by reducing the size of the diagram which defines the homology objects, we obtain a new way to calculate homology. Our main aim is to show Theorem~\ref{Theorem-Endp} which states that the $(n+1)$-st homology of a $k$-extension $A$ may be computed as a limit over the category $\End p$ of all endomorphisms of an $n$-presentation $p$ of $A$.

\begin{notation}\label{Notation-Endp}
 For any $n$-extension $f$ of a $k$-extension $A$, let $\End f$, the \defn{category of endomorphims of $f$ over $A$}, be the full subcategory of $\Ext^{k+n}_A\!\A$ determined by the object $f$. Thus maps in $\End f$ are maps from $f$ to itself which restrict to the identity on $A$ under the functor $\cod^n$.
\end{notation}

When $\A$ has enough projectives we can interpret Proposition~\ref{Proposition-Satellites} the other way round to give

\begin{theorem}[Hopf Formula]\label{Theorem-Hopf-Formula-New}
Let $\A$ be a semi-abelian category with enough projectives and $I\colon {\A\to \B}$ a reflector onto a Birkhoff subcategory $\B$ of $\A$. Let $n\geq 1$. Given an $n$-fold presentation $p$ of an object $A\in|\Extk\A|$, we have
\[
H_{n+1}(A,I_k)\iso K^{n+1}[I_{k+n}p\to T_{k+n}p].
\]
\end{theorem}
\begin{proof}
This is just Proposition~\ref{Proposition-Satellites} viewed from the perspective of Definition~\ref{Definition-Homology}.
\end{proof}

\begin{remark}\label{Remark-Baer-Invariants}
In \cite{Tomasthesis, EverHopf} Everaert gives a direct proof that the right hand side of the Hopf formula is a \textbf{Baer invariant} of $A$: an expression independent of the chosen $n$-fold presentation $p$ of $A$ (see also \cite{EverVdL1, Froehlich}). More precisely, any morphism ${p\to p}$ over $A$ induces the identity on ${K^{n+1}[I_{k+n}p\to T_{k+n}p]}$.
\end{remark}

Of course we can still calculate homology as a limit, as defined in Section~\ref{Section-Without-Projectives}. It turns out that in this case, homology may also be computed as a limit over the small subdiagram of shape $\Init p$, which is a subcategory of $(A\downarrow \cod^n)$. 

\begin{notation}\label{Notation-Initn}
Let $p$ be an $n$-presentation of a $k$-extension $A$. The category $\Init p$ we want to consider is inspired by a higher-dimensional variation on Diagram~\eqref{Diagram-Init}: it is the subcategory of $(A\downarrow \cod^n)$ that is generated by the objects $(p,1_A)$, $(!_P,!_A)$ and $(1_0,!_A)$, all endomorphisms of $p$ over $A$, and the three maps
\[
\vcenter{\xymatrix{P \ar@{ >>}[r]^-p \ar@{=}[d]_{1_P} \ar@{}[dr]|{\Downarrow} & Q \ar[d]^{!_Q} & A \ar[d]_{!_A} \ar@{=}[r]^-{1_A} \ar@{}[dr]|{\Downarrow} & A \ar@{=}[d]^{1_A}\\
P\ar@{ >>}[r]^{!_P} \ar@<.5ex>[d] \ar@{}[dr]|{\Uparrow\Downarrow}& 0 \ar@<.5ex>[d] & 0 \ar@<.5ex>[d] \ar@{}[dr]|{\Uparrow\Downarrow} & A  \ar[l]_-{!_A} \ar@{=}[d] \\
0 \ar@<.5ex>[u] \ar@{=}[r]_-{1_0} & 0 \ar@<.5ex>[u] & 0 \ar@<.5ex>[u] & A \ar[l]^-{!_A}}}
\]
in $(A\downarrow \cod^n)$. The object $A$ is a $k$-extension, but $P$ and $Q$ are  $(k+n-1)$-extensions, with $A$ being the ``terminal object'' of $Q$ (when $Q$ is considered as a diagram in $\Extk\A$). $Q$, the codomain of $p$, is an $(n-1)$-presentation of $A$ (cf.\ definition of $n$-presentation in Subsection~\ref{Subsection-Projective-Presentations}). Note that there is an obvious inclusion ${\End p\to \Init p}$ sending $p$ to $(p,1_A)$.
\end{notation}

\begin{proposition}\label{Proposition-Initnp}
Consider $k\geq 0$, $n\geq 1$ and $A\in|\Extk\A|$, and let $p$ be an $n$-fold presentation of~$A$. Then
\[
K^{n+1}[I_{k+n}p\to T_{k+n}p]=\lim \bigl(\ker^n\comp I_{k+n}\comp U\colon{\Init p\to \Extk\A}\bigr).
\] 
\end{proposition}
\begin{proof}
Kernels and limits commute, so the above limit is also $K^{n-1}[\lim \ker I_{k+n}U]$. Note that $K^{n+1}[I_{k+n}p\to T_{k+n}p]$ is the same as $K^n[I_{k+n}[p]\to T_{k+n}[p]]$, where $I_{k+n}[p]$ and $T_{k+n}[p]$ denote the domains of $I_{k+n}p$ and $T_{k+n}p$ respectively. Thus we only have to show that the limit of $\ker  I_{k+n} U$ coincides with the kernel of $r_p^{k+n}\colon{I_{k+n}[p]\to T_{k+n}[p]}$.

The diagram $\ker I_{k+n}U$ we are considering is
\[
\xymatrix{K[I_{k+n}p]\ar[rrr]^-{\ker(I_{k+n}(1_p,!_Q))=f} &&& K[I_{k+n}!_P]=I_{k+n-1}P \ar@<.5ex>[rrr] &&& 0 \ar@<.5ex>[lll]}
\]
where we name the non-zero map $f$, for convenience. Recall from Remark~\ref{Remark-Bang-Reflector} that $K[I_{k+n}!_P]=I_{k+n-1}P$. We will show that $K[f]=K[r_p^{k+n}]$, which will in turn imply that $K[r_p^{k+n}]$ is indeed the limit of this diagram.

Consider the following diagram, where we are taking kernels to the left. The kernel objects of $T_{k+n}p$ and $I_{k+n-1}p$ are equal, because the bottom right square is a pullback, by definition of $T_{k+n}p$.
\[
\xymatrix@C=8ex{0 \ar[r] & K[I_{k+n}p] \ar[d]_-{\overline{r}}  \ar@{ |>->}[r]^-{\Ker I_{k+n}p} & I_{k+n}[p] \ar@{ >>}[r]^-{I_{k+n}p} \ar@{ >>}[d]_-{r^{k+n}_p} & Q \ar@{=}[d] \ar[r] & 0\\
0 \ar[r] & K[T_{k+n}p] \ar@{ |>->}[r] \ar@{=}[d] & T_{k+n}[p] \ar@{}[rd]|<{\pullback} \ar@{ >>}[r]^-{T_{k+n}p} \ar@{ >>}[d]_-{\overline{\eta}} & Q \ar@{ >>}[d]^{\eta^{k+n-1}_Q} \ar[r] & 0\\
0 \ar[r] & K[I_{k+n-1}p] \ar@{ |>->}[r]_-{\Ker I_{k+n-1}p} & I_{k+n-1}P \ar@{ >>}[r]_{I_{k+n-1}p} & I_{k+n-1}Q \ar[r] & 0}
\]
Recall Diagram~\eqref{Diagram-Compare-I-T}: $\overline{\eta}\comp r_p^{k+n}=\eta^{k+n-1}_{I_{k+n}[p]}$. Looking at the top two exact sequences, we see that the top left square is a pullback, because the arrow $Q\to Q$ at the right hand side is a monomorphism. Thus $K[r^{k+n}_p]=K[\overline{r}]$, and $\overline{r}$ is also a regular epi. We will show that $f$ factors as $f=\Ker I_{k+n-1}p \comp \overline{r}$, and then $K[f]=K[\overline{r}]=K[r^{k+n}_p]$ as desired.

The map $f$ is induced by the following diagram:
\[ 
\xymatrix{K[I_{k+n}p] \ar@{ |>->}[r] \ar[d]^f & I_{k+n}[p] \ar@{ >>}[r]^(.6){I_{k+n}p} \ar@{ >>}[d]^{\eta^{k+n-1}_{I_{k+n}[p]}} & Q \ar[d]\\
I_{k+n-1}P \ar@{=}[r] & I_{k+n}[!_P] \ar[r]_(.6){I_{k+n}!_P} & 0}
\]
It is easily checked that the map $I_{k+n}[p]\to I_{k+n}[!_P]$ above is indeed the same as in the following diagram (c.f.\ also Lemma~\ref{Lemma-II_1=I}):
\[
\xymatrix{K[I_{k+n}p] \ar@{ |>->}[r] \ar@{.. >>}[d]^{\overline{r}} \ar[dr]^f & I_{k+n}[p] \ar@{ >>}[r]^{I_{k+n}p} \ar@{ >>}[d]^{\eta^{k+n-1}_{I_{k+n}[p]}} & Q \ar@{ >>}[d]^{\eta^{k+n-1}_{Q}} \\
K[I_{k+n-1}p] \ar@{ |>->}[r] & I_{k+n-1}P \ar@{ >>}[r]_{I_{k+n-1}p} & I_{k+n-1}Q}
\]
Thus $f$ factors as promised and we have $K[f]=K[r^{k+n}_p]$.

A cone $(C,\sigma)$ on $\ker I_{k+n}U\colon{\Init_{n}p\to \Extk\A}$ consists of three maps $\sigma_{(p,1_A)}$, $\sigma_{(!_P,!_A)}$ and $\sigma_{(1_0,!_A)}$:
\[
\xymatrix@!0@=4,5em{C \ar[d]_-{\sigma_{(1_{0},!_{A})}} \ar[drr]|-{\sigma_{(!_P,!_{A})}} \ar[rr]^-{\sigma_{(p,1_{A})}} && K[I_{k+n}p] \ar[d]^-{f}\\
0 \ar@<.5ex>[rr] && I_{k+n-1}P. \ar@<.5ex>[ll]}
\]
We see that $f\comp \sigma_{(p,1_{A})}=0$, so $\sigma_{(p,1_{A})}$ factors over
\[
\Ker f= \Ker \overline{r}\colon{K[r^{k+n}_{p}]\to K[I_{k+n}p]}.
\]
This factorisation is the needed map ${C\to K[r_p^{k+n}\colon {I_{k+n}[p]\to T_{k+n}[p]}]}$.

We still have to show that $K[r_p^{k+n}]$ itself forms a cone on the diagram $\ker I_{k+n}U$.  Suppose $(g,h)$ is any endomorphism of $p$ over $A$, and write $\overline{g}$ for the induced morphism $\ker(I_{k+n}(g,h))\colon {K[I_{k+n}p]\to K[I_{k+n}p]}$. We have to check that $\overline{g}\comp \Ker\overline{r}=\Ker \overline{r}$. This, however, is a consequence of the fact that $K[r_p^{k+n}]$ is a Baer invariant of~$A$---see Remark~\ref{Remark-Baer-Invariants}. Indeed, in the diagram of short exact sequences
\[
\xymatrix@C=8ex{0 \ar[r] & K[r_p^{k+n}] \ar@{ |>->}[r]^{\Ker r_p^{k+n}} \ar@{.>}[d]_-{g'} & I_{k+n}[p] \ar@{ >>}[r]^-{r^{k+n}_{p}} \ar[d]_-{I_{k+n}[(g,h)]} & T_{k+n}[p] \ar[d]^{T_{k+n}[(g,h)]} \ar[r] & 0\\
0 \ar[r] & K[r_p^{k+n}] \ar@{ |>->}[r]_{\Ker r_p^{k+n}} & I_{k+n}[p] \ar@{ >>}[r]_-{r^{k+n}_{p}} & T_{k+n}[p] \ar[r] & 0}
\]
the induced arrow $g'$ is $1_{K[r_p^{k+n}]}$; hence, using that $\Ker r_p^{k+n}=\Ker I_{k+n}p\comp \Ker \overline{r}$, we get
\[
\Ker I_{k+n}p\comp \Ker \overline{r}=\Ker I_{k+n}p\comp h\comp \Ker\overline{r} 
\]
and the needed equality follows.
\end{proof}

\begin{theorem}\label{Theorem-Endp}
Consider $k\geq 0$, $n\geq 1$ and $A\in|\Extk\A|$. If $\A$ has enough projectives and $p$ is an $n$-fold presentation of $A$ then
\[
H_{n+1}(A,I_k)=\lim \bigl(\ker^n\comp I_{k+n}\colon{\End p\to \Extk\A}\bigr).
\] 
\end{theorem}
\begin{proof}
By Theorem~\ref{Theorem-Hopf-Formula-New} the $(n+1)$-st homology of $A$ is ${K^{n+1}[I_{k+n}p\to T_{k+n}p]}$. Hence by Proposition~\ref{Proposition-Initnp} it suffices to show that $\End p$ is initial in $\Init_{n}p$. We must check that the slice categories $(\End p\downarrow (p,1_A))$, $(\End p\downarrow (!_P,!_A))$ and $(\End p \downarrow(1_0,!_A))$ are non-empty and connected (here we view $\End p$ as the full subcategory of $\Init p$ determined by $(p,1_A)$). There is only one possible map from $(p,1_A)$ to $(1_0,!_A)$, and the other two categories fulfil the needed conditions essentially because $((1_{P},1_{Q}),(1_{A},1_{A}))$ is a terminal object of the slice category $(\End p \downarrow(p,1_{A}))$, and the only maps in $\Init p$ from $(p,1_A)$ to $(!_P,!_A)$ are compositions of an endomorphism of $(p,1_A)$ with $((1_P,!_Q),(!_A,1_A))$.
\end{proof}

\begin{remark}\label{Remark-Fixed-Points}
This means that \emph{computing the homology of an object essentially amounts to finding fixed points of endomorphisms of a projective presentation of this object.} The use of this technique will be illustrated in Examples~\ref{Example-Cyclic-Groups} and~\ref{Example-Generators-and-Relations}.
\end{remark}

\begin{remark}\label{Remark-Homology-and-Baer-Invariants}
We now come back to Remark~\ref{Remark-Baer-Invariants} and interpret Definition~\ref{Definition-Homology} in terms of Baer invariants. It provides an alternative answer to the following question: ``Given a functor $I\colon{\A\to\A}$ and an object $A$ of $\A$, how can we construct an object $H_{n+1}(A,I)$ out of the $n$-extensions of $A$ in a manner which is independent of any particular chosen extension of $A$?'' The classical example is the Hopf formula
\[
H_2(I,A)_{\G}\cong K^{2}[I_{1}p\to T_{1}p]
\]
which expresses $H_2(A,I)_{\G}$ in terms of a projective presentation $p\colon{P\to A}$ of $A$. Of course, the very existence of the isomorphism implies that the expression on its right hand side cannot depend on the choice of $p$. The idea behind Definition~\ref{Definition-Homology} is different but straightforward: simply take the limit of \emph{all} extensions of $A$. The independence might now be understood as follows. If $p$ is an $n$-presentation of $A$ then $H_{n+1}(A,I)$ is the limit of $\ker^n\comp I_n\colon{\End p\to \A}$, which means that $H_{n+1}(A,I)$ is the universal object with the property that all endomorphisms of $p$ are mapped to the same automorphism of this object, its identity.
\end{remark}

Finally we show, as worked out examples, that we can retrieve well-known results in group homology using our new definition.

\begin{example}[Finite cyclic groups]\label{Example-Cyclic-Groups}
We use the methods of our theory to calculate $H_2(C_n,\ab)$ for any $n\in\N$, where $C_n$ is the cyclic group of order $n$. As $\Z$ is projective and abelian, the map $p\colon {\Z\to C_n}$ which sends $1\in \Z$ to a generator $c\in C_n$ is a projective presentation of $C_n$, and central. Thus $H_2(C_n,\ab)$ is the limit of the diagram $\ker \colon {\End p \to \Gp}$. Now any endomorphism of $p$ must be 
\[
\xymatrix{\Z \ar[d]_-{\cdot(nk+1)} \ar@{ >>}[r]^-{p} \ar@{}[dr]|{\Downarrow} & C_n \ar@{=}[d] \\ 
\Z \ar@{ >>}[r]_-p & C_n}
\]
i.e., multiplication by $(nk+1)$ for some $k\in \Z$. So $H_2(C_n,\ab)$ is the limit of the diagram which has as only object $n\Z$, and maps $\cdot(nk+1)\colon {n\Z \to n\Z}$. If $\lambda\colon {H_2(C_n,\ab)\to n\Z}$
is the leg of the limit cone, we must have $\lambda(x)\cdot (nk+1)=\lambda(x)$ for every element $x\in H_2(C_n,\ab)$ and every $k$. \emph{So we are looking for fixed points of the map $\cdot(nk+1)$.} But as, in $n\Z$, $0$ is the only fixed point of multiplication by $(nk+1)$ for all $k\neq 1$, we have $\lambda(x)=0$ for all $x\in H_2(C_n,\ab)$. Thus, as $\lambda$ is a limit cone and so a monomorphism, $H_2(C_n,\ab)=0$. 
\end{example}

\begin{example}[Generators and relations]\label{Example-Generators-and-Relations}
Given a presentation of a group in terms of generators and relations, for example
\[
A=\langle a_1,\ldots,a_n\ |\ r_i=1\rangle
\]
for some relations $r_i$, the kernel of the free presentation
\[
p\colon F_n\to A
\]
is generated by the relations $r_i$ as a \emph{normal} subgroup of $F_n$. Here $F_n$ is the free group on $n$ generators. But when we go to the centralisation 
\[
\centr p\colon {\frac{F_n}{[K[p],F_n]}\to A},
\]
 every element of the kernel commutes with every other element, so now $K[\centr p]$ is generated by the relations $r_i$ as a subgroup of ${F_n}/{[K[p],F_n]}$. Every endomorphism of $p$ over $A$ must send a generator $a_i$ to $a_ik_i$ for some $k_i\in K[f]$, and any choice of $k_i$ gives such an endomorphism. Thus on $\centr p$ we get endomorphisms that send $a_i\in {F_n}/{[K[p],F_n]}$ to $a_i\prod_j r_j^{\alpha_{ij}}$, for some $\alpha_{ij}\in \Z$, and again any choice of $\alpha_{ij}$ gives an endomorphism. Note that $K[\centr p]$ is an abelian group, since it is in the centre of ${F_n}/{[K[p],F_n]}$. From here it is relatively easy to find the fixed points of the induced endomorphism of $K[\centr p]$, given a specific group in terms of generators and relations. We give as an example
\[
C_n\times C_n=\langle a,b\ |\ a^n=1=b^n, aba^{-1}b^{-1}=1\rangle.
\]
Here $p\colon {F_2\to C_n\times C_n}$, and $K[\centr p]$ is generated by $x=a^n$, $y=b^n$ and $z=aba^{-1}b^{-1}$. Note that as $aba^{-1}b^{-1}$ commutes with everything, we get $(aba^{-1}b^{-1})^n=ab^na^{-1}b^{-n}$, and as $b^n$ also commutes with everything, we have $z^n=1$. As described above, any endomorphism of $\centr p$ induced by one on $p$ sends $a\in {F_2}/{[K[p],F_2]}$ to $ax^{\alpha_1}y^{\alpha_2}z^{\alpha_3}$ and $b$ to $bx^{\beta_1}y^{\beta_2}z^{\beta_3}$. On $K[\centr p]$ this gives 
\begin{align*}
x&\mapsto x^{n\alpha_1+1}y^{n\alpha_2}\\
y&\mapsto x^{n\beta_1}y^{n\beta_2+1}\\
z&\mapsto z
\end{align*}
as the $x$, $y$ and $z$ commute with everything, and $z^n=1$. For $x^{l_1}y^{l_2}z^{l_3}$ to be a fixed point for any of these endomorphisms, we need 
\begin{align*}
l_1\alpha_1+l_2\beta_1 &=0\\
l_1\alpha_2+l_2\beta_2 &=0
\end{align*}
for any choice of $\alpha_i$ and $\beta_i$, or in other words we need
\[
l_1\alpha+l_2\beta=0
\]
for any choice of $\alpha$ and $\beta$. Hence $l_1=l_2=0$, and we have fixed points $z^{l_{3}}$. Since $z^n=1$, we get
\[
H_2(C_n\times C_n,\ab)=C_n.
\]

Note that we can use the diagram over $\Init p$ instead of $\End p$ to see that any fixed point must be of the form $aba^{-1}b^{-1}$ for some $a$ and $b$ (or a product of such), since the fixed point must be sent to the identity in $\ab F_n={F_n}/{[F_n,F_n]}$. 
\[
\xymatrix{H_2(A,\ab) \ar[r] \ar[dr] \ar[d] & K[\centr p] \ar[d]\\
0 \ar[r] & \ab F_n}
\]
Comparing this to the Hopf formula
\[
H_2(A,\ab)=\frac{[F_n,F_n]\cap K[p]}{[K[p],F_n]},
\]
we see that the calculation using our method is exactly the same as the one using the Hopf formula; the only thing that is different is the \emph{interpretation} of these elements as fixed points of certain endomorphims. Note that we of course proved in Proposition~\ref{Proposition-Initnp} that the limit of the diagram $\ker\comp I_1\colon {\Init p\to \A}$ is the expression of the Hopf formula, so this is exactly what you would expect. 
\end{example}

\section*{Acknowledgements}

 We would like to thank George Janelidze for pointing us to his work on satellites as an alternative approach to homology and for fruitful discussions. Thanks also to Alexander Frolkin and James Griffin for helping us with the examples, and to Tomas Everaert, Martin Hyland and Peter Johnstone for useful comments and suggestions.

\bibliography{tim}
\bibliographystyle{amsplain}
\bigskip

\end{document}